\documentclass[onefignum,onetabnum,onealgnum,oneeqnum,reqno]{siamart250211}
\usepackage[papersize={8.5in,11in},margin=1.0in,footskip=0.3in]{geometry}
\usepackage[final]{microtype}
\usepackage{amssymb,amsmath}
\usepackage{bm}
\usepackage{enumitem}
\usepackage{xstring}
\usepackage{booktabs}
\usepackage{multirow}
\usepackage{MnSymbol}
\usepackage{orcidlink}

\setlist[enumerate]{leftmargin=2.6em,label=(\roman*),topsep=0.5em,parsep=0.25em}
\setlist[itemize]{leftmargin=2em,topsep=0.25em,parsep=0.25em}

\definecolor{bluey}{HTML}{008fba}
\hypersetup{urlcolor=bluey,linkcolor=bluey,citecolor=bluey}

\usepackage{silence}
\usepackage[font={sf,small},labelsep=period,labelfont=bf,margin=6mm,skip=1mm,belowskip=-1mm]{caption}

\definecolor{bluesection}{HTML}{1e4b5e}
\makeatletter
\renewcommand\section{\@startsection{section}{1}{.25in}{1.3ex \@plus .5ex \@minus .2ex}{-.5em \@plus -.1em}{\reset@font\normalsize\bfseries\color{bluesection}}}
\renewcommand\subsection{\@startsection{subsection}{2}{.25in}{1.3ex\@plus .5ex \@minus .2ex}{-.5em \@plus -.1em}{\reset@font\normalsize\bfseries\color{bluesection}}}
\renewcommand\subsubsection{\@startsection{subsubsection}{3}{.25in}{1.3ex\@plus .5ex \@minus .2ex}{-.5em \@plus -.1em}{\reset@font\normalsize\bfseries\color{bluesection}}}
\makeatother

\usepackage{algorithm}
\usepackage[noend]{algpseudocode}
\makeatletter
\algrenewcommand\ALG@beginalgorithmic{\footnotesize}
\renewcommand{\ALG@name}{\small Algorithm}
\makeatother

\crefname{figure}{Fig.\!}{Figs.\!}
\crefname{section}{\S\!}{\S\!}
\crefname{subsection}{\S\!\!}{\S\!\!}
\newcommand{\R}{\mathbb R}
\DeclareMathOperator{\diag}{diag}
\DeclareMathOperator*{\argmin}{arg\,min}
\newcommand{\vu}{\mathbf u}
\newcommand{\vn}{\mathbf n}
\newcommand{\vf}{\mathbf f}
\newcommand{\vg}{\mathbf g}
\newcommand{\vh}{\mathbf h}
\newcommand{\vsigma}{\bm\sigma}
\newcommand{\trans}{{\mathsf{T}}}
\newcommand{\jump}[1]{[\![ #1 ]\!]}
\newcommand{\polydeg}{{\smash{\wp}}}
\newcommand{\xin}{x}
\newcommand{\xout}{\bar{x}}
\newcommand{\Rey}{\mathsf{Re}}
\newcommand{\zetaI}{\zeta_{\bm{\sigma}}}
\newcommand{\zetaII}{\zeta_{\bm{u}}}
\newcommand{\omegaI}{\omega_{\bm{u}}}
\newcommand{\omegaII}{\omega_p}
\newcommand{\ball}{\raisebox{-0.2ex}{\includegraphics[height=1.8ex]{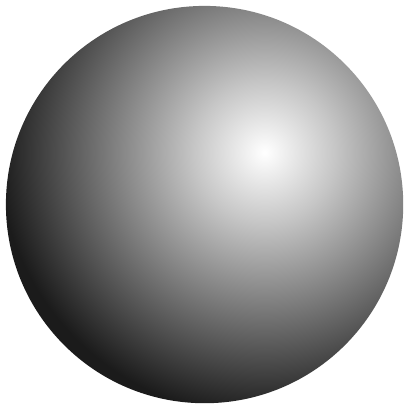}}}
\newcommand{\cube}{\raisebox{-0.65ex}{\includegraphics[height=2.5ex]{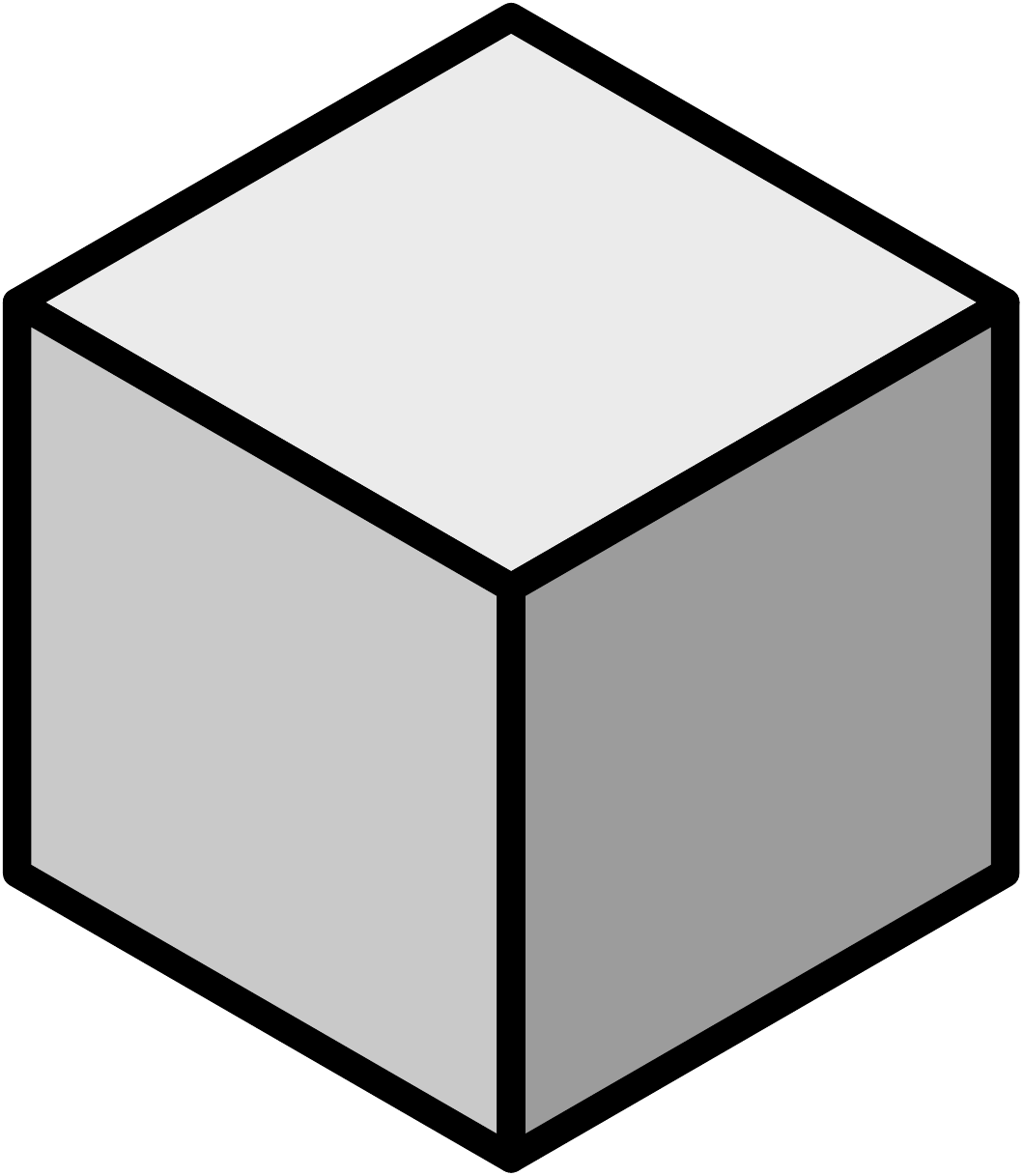}}}
\newcommand{\arcstar}{\raisebox{-0.2ex}{\includegraphics[height=2ex]{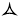}}}
\newcommand{\bcP}{\arcstar}
\newcommand{\bcN}{\raisebox{-0.1ex}{\scalebox{1.2}{$\circ$}}}
\newcommand{\bcD}{\raisebox{-0.1ex}{\scalebox{1.2}{$\bullet$}}}
\newcommand{\Vup}{\mathcal V}
\newcommand{\bg}{\mathring{G}}
\newcommand{\dt}{\delta}

\title{On efficient multigrid solvers for mixed-degree\\local discontinuous Galerkin multiphase Stokes problems}
\author{Robert I.~Saye \orcidlink{0000-0001-7256-6941}}
\date{\today}

\AddToHook{cmd/appendix/before}{\crefalias{section}{appendix}}
\begin{document}
\setlength{\floatsep}{5mm}
\setlength{\textfloatsep}{5mm}

\makeatletter
\begin{center}
{\bfseries\MakeUppercase\@title\par}%
\vskip2mm%
{\small\textsc{\@author}\par}%
\vskip1mm%
{\sffamily\footnotesize Lawrence Berkeley National Laboratory, Berkeley, California, USA $\mid$ \texttt{rsaye@lbl.gov} $\mid$ September 12, 2026 \par}%
\vskip5mm%
\end{center}
\makeatother

\begin{abstract}
We design and investigate efficient multigrid solvers for multiphase Stokes problems discretised via mixed-degree local discontinuous Galerkin methods. Using the template of a standard multigrid V-cycle, we develop a smoother analogous to element-wise block Gauss-Seidel, except the diagonal block inverses are replaced with an approximation that balances the smoothing of the velocity and pressure variables, factoring in the unequal scaling of the various Stokes system operators, and optimised via two-grid local Fourier analysis. We evaluate the performance of the multigrid solver across an extensive range of two- and three-dimensional test problems, including steady-state and unsteady, standard-form and stress-form, single-phase and high-contrast multiphase Stokes problems, with multiple kinds of boundary conditions and various choices of polynomial degree. In the lowest-degree case, i.e., that of piecewise constant pressure fields, we observe reliable multigrid convergence rates, though not especially fast. However, in every other case, we see rapid convergence rates matching those of classical Poisson-style geometric multigrid methods; e.g., 5 iterations reduce the Stokes system residual by 5 to 10 orders of magnitude.
\end{abstract}

\begin{keywords}
multigrid, Stokes, local discontinuous Galerkin, mixed-degree, multiphase, local Fourier analysis
\end{keywords}

\begin{MSCcodes}
65N55 (primary), 65N30, 65N22, 76D07, 65F10
\end{MSCcodes}

\section{Introduction}

In this work, we design and investigate efficient multigrid solvers for multiphase Stokes problems discretised via mixed-degree local discontinuous Galerkin (LDG) methods. Our motivation stems in part from prior work on multigrid Stokes solvers for equal-degree LDG methods, i.e., those for which the velocity and pressure fields have the same polynomial degree \cite{flame,inferno}. In the equal-degree setting, and for various kinds of Stokes problems, one can build highly effective multigrid solvers via simple multigrid V-cycles and element-wise block Gauss-Seidel smoothers; e.g., to reduce the residual by 8 orders of magnitude, only 6 to 8 solver iterations are needed. However, these multigrid solvers crucially depend on a carefully-tuned mechanism of pressure penalty stabilisation, but for some kinds of Stokes problems---notably, those with prescribed stress boundary conditions, or multiphase problems with interfacial stress jump conditions---there is no good way to implement this mechanism (as discussed further in \cref{sec:comparison}).

In contrast, a mixed-degree LDG framework does not require any pressure penalty stabilisation, can solve the aforementioned Stokes problems, and, perhaps counterintuitively, can yield significantly more accurate numerical solutions. %
Here, mixed-degree refers to the piecewise polynomial pressure field being of one degree less than the velocity field. %
Unfortunately, however, the above-mentioned multigrid solvers are ineffective in a mixed-degree setting---{for example, on an elementary Stokes problem, a simple block Gauss-Seidel multigrid solver requires 6 iterations for a $10^8$ residual reduction in the equal-degree setting, but fails in the mixed-degree setting unless the smoother is damped; even with optimal damping, convergence is slow and degrades as the mesh is refined (see \cref{app:briefcomparison}).} Consequently, the main aim of this study is to develop a replacement multigrid solver, specifically by targeting the multigrid smoother.

One possibility is to apply sparse approximate inverses (SAI) methods, as explored in \cite{inferno}. These techniques have a number of potential benefits, including high parallelism, increased robustness (e.g., in handling a wider range of discretisation parameters), and the potential for black-box adaptivity (e.g., to handle anisotropy, highly unstructured meshes, and complex geometry, etc.). SAI approaches can lead to rapid multigrid convergence; on the other hand, they can be relatively expensive to build and apply. In this work, we leverage some of the techniques developed in \cite{inferno}, but with the aim of devising a smoother whose application cost is the same as block Gauss-Seidel. In essence, the developed smoother is identical to multicoloured block Gauss-Seidel, except the diagonal block inverses are replaced with an approximation that optimally balances the various Stokes system operators (i.e., the viscous, pressure gradient, velocity divergence, and density-weighted temporal derivative operators). Given its similarity to element-wise block Gauss-Seidel, one may also view the smoother as a multiplicative Schwarz method, whereby the subdomains of the Schwarz method correspond to the combined set of degrees of freedom on each mesh element; likewise, the smoother shares various operational similarities with that of multiplicative Vanka-type methods \cite{Vanka1986,Vanka1986b,doi.org/10.1002/nla.2306} and box relaxation methods \cite{BrandtLivne}, though differs in how the smoother is constructed.

The paper is structured as follows. In \cref{sec:ldg}, we outline the prototype Stokes problem and summarise the mixed-degree LDG framework; we also include numerical examples and convergence tests to compare against the equal-degree setting. In \cref{sec:mg}, we develop the multigrid framework, starting with basic ingredients of the V-cycle and then focusing on the smoother, including: (i) its construction, (ii) balancing methods to treat the unequal scaling of the Stokes system operators, (iii) adaptation to unsteady Stokes problems, (iv) overall assembly and application of the smoother, and (v) optimisation of its parameters. Numerical results are then presented in \cref{sec:experiment}, examining multigrid performance across a variety of test problems, including 2D and 3D steady-state and unsteady Stokes problems in standard-form and stress-form, different kinds of boundary conditions, and high-contrast multiphase problems. A concluding perspective is given in \cref{sec:conclusion}, while additional details of the LDG framework and smoother parameter optimisation are given in \cref{app:ldg,app:lfa}, respectively.

\section{LDG Scheme}
\label{sec:ldg}

Our prototype (steady-state) Stokes problem consists of solving for a velocity field $\vu : \Omega \to \R^d$ and pressure field $p: \Omega \to \R$ such that
\begin{equation} \label{eq:govern1} \left. \begin{aligned} - \nabla \cdot \bigl(\mu_i (\nabla \vu + \gamma\,\nabla \vu^\trans) \bigr) + \nabla p &= \vf \\ -\nabla \cdot \mathbf \vu &= f \end{aligned} \right\} \text{ in } \Omega_i, \end{equation}
subject to the interfacial jump conditions
\begin{equation} \label{eq:govern2} \left. \begin{aligned} \jump{\vu} &= \vg_{ij} \\ \jump{\mu (\nabla \vu + \gamma\,\nabla \vu^\trans) \cdot \vn - p \vn} &= \vh_{ij} \end{aligned} \right\} \text{ on } \Gamma_{ij}, \end{equation}
and boundary conditions
\begin{equation} \label{eq:govern3} \begin{aligned} \vu &= \vg_{\partial} && \text{on } \Gamma_D, \\ \mu (\nabla \vu + \gamma\,\nabla \vu^\trans) \cdot \vn - p \vn &= \vh_{\partial} && \text{on } \Gamma_N, \end{aligned} \end{equation}
where $\Omega$ is a domain in $\mathbb R^d$ divided into one or more subdomains/phases $\Omega_i$, $\mu_i > 0$ is a phase-dependent viscosity coefficient, $\Gamma_{ij} := \partial \Omega_i \cap \partial \Omega_j$ is the interface between phase $i$ and $j$, and $\Gamma_D$ and $\Gamma_N$ denote the parts of $\partial \Omega$ upon which velocity Dirichlet or stress boundary conditions are prescribed, respectively. Here, $\gamma \in \{0,1\}$ is a parameter defining the form of the Stokes equations: they are said to be in \textit{standard form} when $\gamma = 0$ and \textit{stress form} when $\gamma = 1$. Meanwhile, $\vf$, $f$, $\vg$, and $\vh$ provide the data for the multiphase Stokes problem and are given functions defined on $\Omega$, its boundary, and internal interfaces.

We consider here meshes arising from Cartesian grids, quadtrees or octrees; accordingly, a tensor-product polynomial space is adopted. Let $\mathcal E = \bigcup_i E_i$ denote the mesh and its elements, let $\polydeg \geq 0$ be an integer, and let $\mathcal Q_\polydeg$ be the space of tensor-product polynomials of one-dimensional degree $\polydeg$; e.g., $\mathcal Q_2$ is the space of biquadratic or triquadratic polynomials in 2D or 3D, respectively. Define the piecewise polynomial space
\[ V_{\polydeg} := \bigl\{ u : \Omega \to \R \enskip \bigl| \enskip u|_E \in \mathcal Q_{\polydeg} \text{ for every } E \in \mathcal E \bigr\}. \]
Note that an arbitrary function of $V_\polydeg$ has $(\polydeg + 1)^d$ degrees of freedom per element. We similarly define $V_\polydeg^d$ and $V_\polydeg^{d \times d}$ to be the space of piecewise polynomial vector-valued and matrix-valued fields, respectively.

The LDG framework considered here solves \cref{eq:govern1,eq:govern2,eq:govern3} for a discrete solution $\bm u_h \in V_\polydeg^d$ and $p_h \in V_{\polydeg-1}$ according to a five-step process: (i) using viscosity-upwinded one-sided numerical fluxes, define a discrete gradient $\bm \tau_h \in V_\polydeg^{d \times d}$ equal to the discretisation of $\nabla \vu_h$, taking into account Dirichlet source data $\vg$; (ii) using $\bm \tau_h$, define a viscous stress tensor $\vsigma_h \in V_\polydeg^{d \times d}$ equal to the $L^2$ projection of $\mu (\bm \tau_h + \gamma\,\bm \tau_h^\trans) - p_h \mathbb I$ onto $V_\polydeg^{d \times d}$; (iii) using the adjoint of the discrete gradient operator from the first step, compute the negative discrete divergence of $\vsigma_h$, taking into account Neumann-like data $\vh$, and set the result equal to the $L^2$ projection of $\vf$ onto $V_\polydeg^d$; (iv) compute a mixed-degree negative discrete divergence of $\vu_h$, taking into account Dirichlet data $\vg \cdot \vn$, and set the result equal to the $L^2$ projection of $f$ onto $V_{\polydeg-1}$; last, (v) amend the overall system by appropriately including penalty stabilisation terms. Our construction is partly based on the LDG schemes initially developed by Cockburn \textit{et al} \cite{CockburnKanschatSchotzauSchwab2002}, with extensions to multiphase variable-viscosity problems analogous to the techniques developed in prior work \cite{flame}; details are deferred to  \cref{app:ldg}. The outcome is a symmetric linear system of the following form (here illustrated for the 3D case, with that of the 2D case being an obvious reduction),
\begin{equation} \label{eq:components}
\begin{pmatrix} \mathcal{A}_{11} & \mathcal{A}_{12} & \mathcal{A}_{13} & -\widetilde{G}_1^\trans \bar M \\ \mathcal{A}_{21} & \mathcal{A}_{22} & \mathcal{A}_{23} & -\widetilde{G}_2^\trans \bar M \\ \mathcal{A}_{31} & \mathcal{A}_{32} & \mathcal{A}_{33} & -\widetilde{G}_3^\trans \bar M \\
-\bar M \widetilde{G}_1 & -\bar M \widetilde{G}_2 & -\bar M \widetilde{G}_3 & 0 \end{pmatrix} \begin{pmatrix} u_h \\ v_h \\ w_h \\ p_h \end{pmatrix} = \begin{pmatrix} b_u \\ b_v \\ b_w \\ b_p \end{pmatrix}\!, \end{equation}
where $u_h, v_h, w_h \in V_\polydeg$ denote the components of $\vu_h$, while $(b_u,b_v,b_w,b_p)$ collects the entire influence of the source data $\vf$, $f$, $\vg$, and $\vh$ onto the right-hand side.\footnote{Scattered throughout this paper we often apply a few convenient abuses of notation. A piecewise polynomial function (e.g., in $V_\polydeg$) may carry the same notation as its corresponding coefficient vector relative to the chosen basis of $V_\polydeg$, with the precise meaning understood from context. For example, in the identity $u^\trans M v = \int_\Omega u\,v$, the left-hand side uses vectors and matrices, whereas the right-hand side applies the functional form. Likewise, operator notation may refer to the functional form or to its matrix representation relative to the chosen basis. %
For example, the $i$th component of the broken gradient operator $\bg : V_\polydeg \to V_\polydeg^d$ satisfies $u^\trans\!M \bg_i v = \int_\Omega u\,\bg_i v = \sum_{E \in \mathcal E} \int_{\!E} u\,\partial_i v$, for every $u,v \in V_\polydeg$.} Here, $\widetilde{G} : V_{\polydeg} \to V_{\polydeg - 1}^d$ is a mixed-degree discrete gradient operator, $\bar M$ is a block-diagonal mass matrix relative to $V_{\polydeg - 1}$, and $\mathcal{A} : V_\polydeg^d \to V_\polydeg^d$ implements the viscous part of the Stokes momentum equations, whose $(i,j)$th block is given by 
\[ \mathcal{A}_{ij} = \delta_{ij} \bigl( \textstyle{\sum_{k=1}^d} G_k^\trans M_\mu G_k \bigr) + \gamma\, G_j^\trans M_\mu G_i + \delta_{ij} \mathcal{E}, \]
where $\delta_{ij}$ is the Kronecker delta, $G: V_\polydeg \to V_\polydeg^d$ is a full-degree discrete gradient operator, $M_\mu$ is a $\mu$-weighted block-diagonal mass matrix relative to $V_\polydeg$, and $\mathcal{E}$ is a block-sparse matrix implementing velocity penalty stabilisation. \Cref{eq:components} gives the system in component form; we will also refer to it through an abridged notation,
\begin{equation} \label{eq:blockform} \begin{pmatrix} \mathcal{A} & \mathcal{G} \\ \mathcal{D} & 0 \end{pmatrix} \begin{pmatrix} \vu_h \\ p_h \end{pmatrix} = \begin{pmatrix} {\mathbf b}_{\vu} \\ b_p \end{pmatrix}\!, \end{equation}
where $\mathcal G : V_{\polydeg-1} \to V_\polydeg^d$ is an effective pressure gradient operator and $\mathcal D : V_\polydeg^d \to V_{\polydeg-1}$ is an effective negated velocity divergence operator. The row $\begin{pmatrix} \mathcal A & \mathcal G \end{pmatrix}$ is referred to as the stress divergence operator and the row $\begin{pmatrix} \mathcal D & 0 \end{pmatrix}$ as the divergence constraint operator. In terms of multigrid design, we shall also refer to the Stokes system \cref{eq:components} and \cref{eq:blockform} through the condensed form
\[ A_h x_h = b_h, \]
where $x_h \in \mathcal{V}_h$ collects $\vu_h$ and $p_h$ into one set of unknowns, $\mathcal{V}_h := V_\polydeg^d \otimes V_{\polydeg-1}$, and $A_h$ is the entire symmetric saddle point operator. (The subscript $h$ may be dropped when context permits.)

\subsection{Comparison to equal-degree framework}
\label{sec:comparison}

Before we shift focus to the development of multigrid solvers, we briefly compare the mixed-degree and equal-degree LDG approaches. The equal-degree framework is essentially the same except for three main aspects: (i) the pressure field is of full degree, i.e., $p_h \in V_\polydeg$; (ii) the operators $\mathcal G$ and $\mathcal D$ in \cref{eq:blockform} are replaced by analogous full-degree counterparts; and (iii) a pressure penalty stabilisation operator $\mathcal{P}$ is included in the bottom-right block of the Stokes system, so that it takes the form $\bigl(\begin{smallmatrix} \mathcal A & \mathcal G \\ \mathcal D & -\mathcal{P} \end{smallmatrix}\bigr)$. 
The purpose of $\mathcal{P}$ is to ensure stability in the discretised Stokes system by weakly enforcing continuity of the pressure field; for equal-degree LDG discretisations it is defined such that \cite{CockburnKanschatSchotzauSchwab2002}
\begin{equation} \label{eq:ppenalty} q^\trans \mathcal{P} p = \int_{\Gamma_\circ} \tau_p\, \jump{q} \jump{p}, \end{equation}
for every $p, q \in V_\polydeg$, where $\Gamma_\circ$ denotes the set of all interior mesh faces and $\tau_p$ is a pressure penalty parameter (proportional to the element size $h$ and inversely proportional to the local viscosity).\footnote{One may also interpret $\mathcal P$ as penalising the difference between $\mathcal{G} p$ and the broken gradient of $p$.} This form of pressure penalisation is applicable whenever the pressure field is expected to be continuous, or, via a simple modification, whenever interfacial jumps in $p$ can be specified \textit{prior} to solving the Stokes problem. However, for multiphase problems involving interfacial stress jump conditions, the interfacial jump in $p$ is not known a priori, because it cannot be inferred from the given stress jump conditions $\jump{\vsigma \cdot \vn}$ without first knowing $\jump{\mu (\nabla \vu + \nabla \vu^\trans) \cdot \vn}$. Similarly, for problems with prescribed stress boundary conditions, we cannot weakly enforce pressure boundary conditions via \cref{eq:ppenalty} (were it to include boundary faces) because $p|_{\partial \Omega}$ is not known a priori. In general, some kind of pressure stabilisation is needed for equal-degree LDG schemes in order to ensure stability of the discretised Stokes system \cite{CockburnKanschatSchotzauSchwab2002}. As part of this work, some alternative pressure penalisation approaches were investigated but none were found satisfactory,\footnote{We mention here two approaches. (i) One possibility is to include a discrete viscous stress tensor $\vsigma_h$ as an additional solution variable, thereby expanding the coupled Stokes system to solve for velocity (a vector field), pressure (a scalar field), and stress (a matrix-valued field);
then, instead of penalising $\jump{p_h}$, one may penalise $\jump{\vsigma_h \cdot \vn}$. This approach results in a stable discretisation and brief experimentation suggests that one can build fast-converging multigrid solvers. However, the substantial increase in system size renders this approach too computationally expensive, especially for high-order 3D problems. (ii) Another approach is to generalise \cref{eq:ppenalty} so as to infer $\jump{p}$ as part of the solution process itself. Two possibilities that preserve the overall symmetry of the Stokes operator are to replace the bilinear penalty functional with either
\[ \int_{\Gamma_\circ} \jump{\vsigma(\mathbf v,q) \cdot \vn} \cdot \jump{\vsigma(\vu,p) \cdot \vn - \vh} \qquad \text{or} \qquad  \int_{\Gamma_\circ} \jump{\vn \cdot \vsigma(\mathbf v,q) \cdot \vn} \jump{\vn \cdot \vsigma(\vu,p) \cdot \vn - \vh\cdot\vn},\]
where $\mathbf v$ and $q$ are test functions and $\vsigma(\cdot,\cdot)$ is an LDG viscous stress tensor operator. In either case, the new Stokes system takes the form
\[ \begin{pmatrix} \mathcal A - \mathcal{P}_{11} & \mathcal G + \mathcal{P}_{12} \\ \mathcal D + \mathcal{P}_{12}^\trans & -\mathcal{P}_{22} \end{pmatrix}, \]
where both $\mathcal{P}_{11}$ and $\mathcal{P}_{22}$ are symmetric positive semidefinite. Unfortunately, neither form of stress penalisation is effective in practice---to stabilise pressure, the penalisation needs to be sufficiently strong; on the other hand, it needs to be sufficiently weak so as to preserve enough of the smoothing properties of $\mathcal A$ needed for effective multigrid solvers; in fact, if it is too strong then $\mathcal A - \mathcal{P}_{11}$ becomes indefinite and stability is lost entirely---numerical experiments indicate there is no good middle ground between these competing needs.} hence the adoption here of a mixed-degree framework that requires no pressure penalisation whatsoever.

\begin{figure}%
\centering%
\includegraphics[scale=1.05]{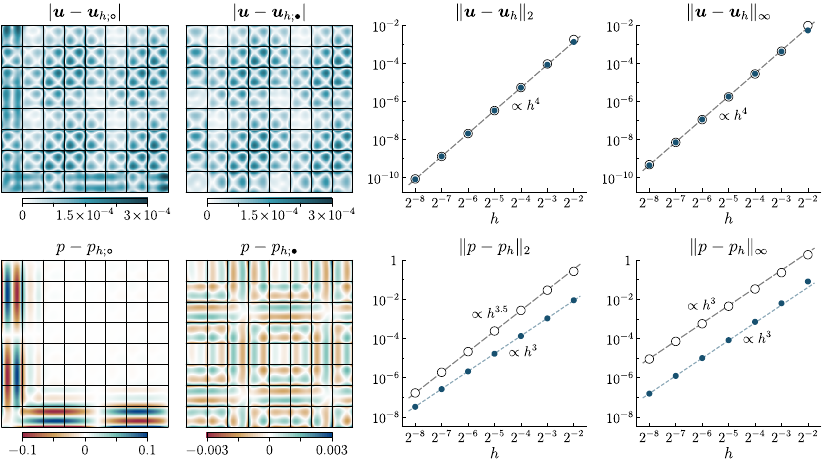}%
\caption{Comparing solution accuracy for the two main LDG approaches: equal-degree (far-left and white circles $\circ$) and mixed-degree (middle-left and filled circles $\bullet$).}%
\label{fig:accuracy}%
\end{figure}

Besides the issue of pressure penalisation, another key difference between the mixed-degree and equal-degree frameworks is that of solution accuracy. In particular, the equal-degree approach often leads to significant pressure errors near boundaries and interfaces, whereas the mixed-degree approach yields cleaner and more accurate pressure fields. To demonstrate this, the left portion of \cref{fig:accuracy} compares discrete solutions arising from a single-phase stress-form Stokes problem, with $\mu \equiv 1$ on $\Omega = (0,1)^2$, stress boundary conditions, and source data generated from a sinusoidal exact solution.\footnote{The exact solution used to generate the data of \cref{fig:accuracy} consists of unit-amplitude sine waves; specifically $\vu = (u,v)$ where $u = u(x,y) = \cos( 2 \pi x + 0.25) \cos(2\pi y + 1.25)$, $v = v(x,y) = \cos(2 \pi x + 4.25) \cos(2 \pi y + 5.25)$, and $p = p(x,y) = \cos(2 \pi x - 3.75) \cos(2 \pi y -2.75)$.} The far-left (resp., middle-left) figures plot the solution error of the equal-degree (resp., mixed-degree) solution, as computed via $\polydeg = 3$ on an $8 \times 8$ mesh. Note that in the equal-degree case, the pressure field contains significant boundary errors, around two orders of magnitude greater than the interior; although the differences in velocity field error are less pronounced, one may also observe a weak numerical boundary layer in the top-left plot. It is important to note that similar outcomes occur for Dirichlet boundary conditions and for embedded interfaces in multiphase problems---as such, the large numerical boundary layers of the equal-degree approach are ingrained in the discretisation itself (and not so much the nuances of penalisation). Using the same test problem, the right-half of \cref{fig:accuracy} contains the results of a grid convergence test, examining the maximum and $L^2$ error norms over grid sizes $4 \times 4$ up to $256 \times 256$. Note that both approaches yield nearly identical velocity errors of optimal order accuracy, i.e., order $\polydeg + 1$. Regarding the pressure error, once again we see a substantial difference---both approaches yield order $\polydeg$ in the maximum norm (for mixed-degree, this is optimal; for equal-degree, this is suboptimal), but the mixed-degree pressure error is smaller by a factor of $\approx 60$. Comparison in the $L^2$ norm is more subtle: roughly speaking, it is possible for the equal-degree setting to yield optimal order pressure fields sufficiently far from boundaries/interfaces, which means it is possible for $\|p - p_{h;\circ}\|_2 \in \mathcal{O}(h^{\polydeg + 1/2})$, as observed in the bottom-middle-right plot of \cref{fig:accuracy}. Nevertheless, the mixed-degree case yields better $L^2$-norm pressure errors (of order $\polydeg$). Overall, to various extents the same observations hold for other kinds of test problems and other choices of polynomial degrees $\polydeg$, making the mixed-degree approach a compelling choice for LDG-based Stokes solvers.

\subsection{Choice of basis} One may implement the above LDG framework using any reasonable basis; for example, a tensor-product Gauss-Lobatto nodal basis or a tensor-product Legendre polynomial orthogonal basis. We have used the latter approach in this work as it simplifies the design and implementation of the multigrid smoother.

\section{Multigrid Schemes}
\label{sec:mg}

The multigrid framework considered here applies the standard concepts of a V-cycle operating on a nested mesh hierarchy. To summarise the key ingredients:
\begin{itemize}
    \item \textit{Mesh hierarchy.} In this work we use quadtrees and octrees to define the finest mesh; the corresponding tree structure naturally defines a hierarchical process by which to agglomerate mesh elements into a nested mesh hierarchy $\mathcal E_h, \mathcal E_{2h}, \mathcal E_{4h}, \ldots$. 
    \item \textit{Interpolation operator.} The interpolation operator $I_{2h}^h : \Vup_{2h} \to \Vup_h$ transfers coarse mesh corrections to a fine mesh and is defined naturally via polynomial injection: $(I_{2h}^h u)|_{E_f} := u|_{E_c}$, where $E_f$ is a fine mesh element and $E_c \supseteq E_f$ is its corresponding coarse mesh element. Note that $I_{2h}^h$ independently acts on each component of the state: on velocity (resp., pressure) components it is a degree $\polydeg$ (resp., $\polydeg - 1$) interpolation operator.
    \item \textit{Restriction operator.} The restriction operator $R_h^{2h} : \Vup_h \to \Vup_{2h}$ transfers a fine mesh residual onto the coarse mesh and is defined via $L^2$ Galerkin projection onto $\Vup_{2h}$. Equivalently, it can be defined via the adjoint of the interpolation operator, $R_h^{2h} := {\mathcal M}_{2h}^{-1} (I_{2h}^h)^\trans {\mathcal M}_h$, where $\mathcal{M}_h$ is the block-diagonal mass matrix relative to $\Vup_h$. 
\end{itemize}
To define the primary Stokes operator $A_h$ on every level of the hierarchy, there are two possible approaches: (i) one may use a ``pure geometric multigrid'' method in which every level is explicitly meshed (e.g., listing out the elements and faces, their connectivity and quadrature, etc.) upon which $A_{2h}, A_{4h}, \ldots$ is built via the LDG discretisation; alternatively, (ii) one may use an operator-coarsening approach that recursively and conveniently automates this process and guarantees the coarse grid problems are discretised in a way that harmoniously matches the fine grid's discretisation. We have used the latter approach in this work---for more details, see \cite{dgmg}, which introduces LDG operator-coarsening methods for Poisson problems, and \cite{flame}, for their extension to multiphase Stokes problems.

\begin{figure}[!t]%
\centering%
\begin{minipage}{0.52\textwidth}%
\begin{algorithm}[H]
	\caption{\sffamily\small Multigrid V-cycle $V({\mathcal E}_h, x_h, b_h)$ with $\nu_1$ pre- and $\nu_2$ post-smoothing steps on mesh ${\mathcal E}_h$ of the hierarchy.}
	\begin{algorithmic}[1]
		\If{${\mathcal E}_h$ is the bottom level}
			\State Solve $A_h x_h = b_h$ using the bottom level direct solver.
		\Else
			\State Apply pre-smoother $\nu_1$ times.
			\State Compute restricted residual, $r_{2h} := (I_{2h}^h)^\trans (b_h - {A}_h x_h)$. \label{line:restrict}
			\State Solve coarse grid problem, $x_{2h} := V({\mathcal E}_{2h}, 0, r_{2h})$.
			\State Interpolate and correct, $x_h \leftarrow x_h + I_{2h}^h x_{2h}$.
			\State Apply post-smoother $\nu_2$ times.
		\EndIf
		\State{\textbf{return} $x_h$.}
	\end{algorithmic}%
	\label{algo:vcycle}%
\end{algorithm}%
\end{minipage}%
\end{figure}

We consider here the template of a standard V-cycle, as shown in \cref{algo:vcycle}.\footnote{\setlength\lineskiplimit{-\maxdimen}\setlength\baselineskip{1.1\baselineskip}Note that $(I_{2h}^h)^\trans$ appears on line \ref{line:restrict} rather than the operator $R_h^{2h}$; this follows from the Galerkin approach, namely that $A_h x_h = b_h$ implements the PDE problem in the ``mass basis'', i.e., the discretised operator, mapping functions from $V_h$ to $V_h$, is given by $M_h^{-1} A_h$. As a piecewise polynomial function in $V_h$, the residual of the fine mesh problem is given by $M_h^{-1} (b_h - A_h x_h)$. In the multigrid V-cycle, this residual is multiplied by $R_h^{2h}$ to define the right-hand side data for the coarse mesh problem: $M_{2h}^{-1} A_{2h} x_{2h} = R_h^{2h} M_{h}^{-1} (b_h - A_h x_h)$, i.e., $A_{2h} x_{2h} = (I_{2h}^h)^\trans (b_h - A_h x_h) =: r_{2h}$.} The design of the smoother is one of our main focuses in this work. Besides the smoother, we have made the following implementation choices:
\begin{itemize}
    \item \textit{Bottom solver.} On the bottom level of the hierarchy (often consisting of just one element spanning the entire domain), a symmetric eigendecomposition of $A$ is precomputed and used as a direct solver. For Stokes systems, oftentimes there is a nontrivial kernel matching that of the non-discretised continuum problem, e.g., constant velocity and/or pressure fields. In these cases, the bottom solver applies a standard pseudoinverse approach to appropriately handle kernel eigenmodes.
    \item \textit{Pre- and post-smoother iteration counts.} There is a balance between choosing $\nu$ relatively high, thus needing fewer (but more expensive) V-cycles, versus choosing $\nu$ relatively low, thus needing more (but cheaper) V-cycles. Prior work in fast LDG multigrid methods \cite{dgmg,fluxx,flame,inferno} has consistently indicated that $\nu_1 \approx 3$ and $\nu_2 \approx 3$ generally yields the best results in terms of fastest solver time to reduce the residual by a fixed factor. Consequently, to keep the focus on key results, we have fixed $\nu_1 = \nu_2 = 3$. %
    \item \textit{Solver acceleration.} Although multigrid can be used as a standalone iterative method, its convergence can be accelerated by using it as a preconditioner in a Krylov method such as Generalised Minimal Residual (GMRES) \cite{demmel,saad}. For the numerical study of this work, we apply multigrid-preconditioned GMRES as the main solver driver, using a single application of the V-cycle of \cref{algo:vcycle} as left-preconditioner.
\end{itemize}

In the remainder of this section we discuss the design of the multigrid smoother, adaptation to unsteady Stokes problems, application in the multigrid V-cycle, and finally the optimisation of its parameters.

\subsection{Smoother preliminaries}
\label{sec:smoother1}

The multigrid smoother developed in this work is in many ways analogous to a multicoloured block Gauss-Seidel smoother, where each \textit{block} corresponds to the collective set of degrees of freedom on each mesh element, i.e., velocity and pressure combined. Accordingly, partition $Ax = b$ so that $\sum_j A_{ij} x_j = b_i$, where $(\cdot)_i$ denotes the set of coefficients on element $i$ and $A_{ij} \in \mathbb{R}^{n \times n}$ denotes the corresponding $(i,j)$th block of $A$, each of size $n := d(\polydeg + 1)^d + \polydeg^d$. Given an input approximation $\xin \in \mathcal{V}_h$, one iteration of block Gauss-Seidel corresponds to sweeping over the elements in some particular order so as to compute an updated approximation $\xout \in \mathcal{V}_h$ whereby
\begin{equation} \label{eq:gssweep} \xout_i := \xin_i + A_{ii}^{-1} \bigl(b_i - \textstyle \sum_{j} A_{ij} \mathring{x}_{j,i} \bigr) \quad \text{where} \quad \mathring{x}_{j,i} := \begin{cases} \xout_j & \text{if } j < i, \\ \xin_j & \text{otherwise,} \end{cases} \end{equation}
where ``$j < i$'' means block $j$ is updated before block $i$. As mentioned in the introduction, for equal-degree LDG Stokes problems, block Gauss-Seidel leads to highly effective multigrid solvers (see, e.g., \cite{flame,inferno}); however, it is ineffective on mixed-degree LDG problems, even with the best possible damping/overrelaxation parameter---V-cycle iteration counts are significantly higher and convergence rates strongly decay as the grid is refined. One of the goals in this work is to develop a replacement smoother that has the same application cost as block Gauss-Seidel; in essence, this is achieved by replacing the inverse block-diagonal term $A_{ii}^{-1}$ in \eqref{eq:gssweep} by a regularised version, denoted $Q_i$, partly inspired by balanced SAI smoothers \cite{inferno}. In turn, $Q_i$ is computed through a weighted least squares approximation so that $A \diag(Q_i) \approx \mathbb I$, except relatively few off-diagonal components of $A$ enter the least squares calculation, the right-hand side is altered to include damping/overrelaxation coefficients, and the weighting is chosen following a balancing process designed to smooth the velocity and pressure variables at the same rate. These aspects are controlled by four positive parameters: $\zetaI$, $\zetaII$, $\omegaI$, and $\omegaII$. The parameters $\zetaI$ (resp., $\zetaII$) control the weighting of the rows corresponding to the stress divergence operator (resp., divergence constraint operator) on the diagonal block of $A$. Meanwhile, $\omegaI$ (resp., $\omegaII$) effectuate the damping or overrelaxation of the velocity (resp., pressure) degrees of freedom, thereby allowing these variables to be damped/overrelaxed by different amounts. The value of these parameters are tunable but depend only on $\polydeg$, $\gamma$, and $d$ and are otherwise constant; typically, $\omegaI$ and $\omegaII$ are close to unity, while $\zetaII$ is unit-order in magnitude, and we shall see later that $\zetaI$ can in essence be eliminated from consideration.

To detail the construction of the smoother $Q$, some additional notation is useful. For an arbitrary coefficient vector $x \in \mathcal{V}_h = V_\polydeg^d \otimes V_{\polydeg-1}$, let $x_{i,\vu,\ell,k}$ denote the component of $x$ associated with the $k$th modal coefficient of the $\ell$th component of the velocity field on the $i$th element; similarly, let $x_{i,p,k}$ denote the $k$th degree of freedom of the pressure field on the $i$th element.\footnote{Beyond the grouping into element-wise blocks and the consistent matching of the degrees of freedom of blockwise matrices and vectors, none of the algorithms described in this work depend on the particular ordering thereof.} In particular, let $k = 0$ correspond to the constant-mode lowest-degree modal coefficient; for example, if $x \equiv 0$ except for $x_{i,p,0} \neq 0$, then $x$ corresponds to an identically zero velocity field and piecewise constant pressure field.

For the moment, suppose that we have appropriately ``balanced'' the Stokes operator. Then, corresponding to the construction of $Q_i$ on the $i$th element, define the weighting vector $\zeta_i \in \mathcal{V}_h$ such that $\zeta_i \equiv 0$, except for the following components: $\zeta_{i;i,\vu,\ell,k} = \zetaI$ for every $\ell$ and $k$; $\zeta_{i;i,p,k} = \zetaII$ for every $k$; and $\zeta_{i;j,\vu,\ell,0} = 1$ for every $\ell$ and $j \neq i$. In addition, define the damping/overrelaxation vector $\omega \in \mathcal{V}_h$ such that $\omega \equiv 0$, except for the following components: $\omega_{j,\vu,\ell,k} = \omegaI$ for every $j,\ell,k$; and $\omega_{j,p,k} = \omegaII$ for every $j,k$. Then, $Q_i$ is defined as the solution of the following minimisation problem:
\begin{equation} \label{eq:Qi} Q_i = \argmin_{\mathsf{Q} \in \mathbb{R}^{n \times n}}~\bigl\| \diag(\zeta_i) \bigl( A_{:,i}\,\mathsf{Q} - \diag(\omega)\, \mathbb{I}_{:,i} \bigr) \bigr\|_F^2, \end{equation}
where $A_{:,i}$ denotes the $i$th block column of $A$ and $\| \cdot \|_F$ denotes the Frobenius norm. \Cref{eq:Qi} implements an overdetermined weighted least squares problem, whose matrix is essentially a trimmed version of the $i$th block column of $A$, with right-hand side data equal to the $i$th block column of the $\omega$-weighted identity matrix $\mathbb{I}_\omega := \diag(\omega) \mathbb I$. Importantly, note that among the nonzero off-diagonal blocks of $A_{:,i}$, only a small subset of rows are actually kept: the weighting vector $\zeta_i$, which multiplies the rows of $A_{:,i}$, is zero everywhere except for the rows of the stress divergence operator associated with constant-mode lowest-degree modal coefficients. There are only $d$ such rows per off-diagonal block and all have unit weighting, while the weighting of the rows in the diagonal block are controlled by $\zetaI$ and $\zetaII$. %

\subsection{Balancing}
\label{sec:balance}

We next consider the balancing of the Stokes operator, via a simple diagonal pre- and post-scaling of $A$. This balancing approach is similar to that developed in prior work \cite{inferno}, and treats two key aspects of Stokes problems:
\begin{itemize}
    \item \textit{Mesh-dependent operator scaling.} Note that the operators within \cref{eq:blockform} scale with $h$ and $\mu$ in different ways: $\mathcal A \sim \mu h^{d-2}$, $\mathcal G \sim h^{d-1}$, and $\mathcal D \sim h^{d-1}$. In particular, note the dependence on $h$ means that the length scales of the domain, or even the particular level of the multigrid hierarchy, can change the relative strengths of these operators, even when the entire collection of elements have the same viscosity coefficient. Without treatment, this may skew the solution of the least squares problem in \cref{eq:Qi}, by over- or under-penalising the off-diagonal components, or by counteracting the role of the weighting parameters $\zetaI$ and $\zetaII$ that are meant to (universally, for any $\mu$ and any $h$) balance the smoothing characteristics of velocity and pressure variables.
    \item \textit{High-contrast problems.} Whenever the viscosity coefficient $\mu$ exhibits large variations or discontinuities, the various (sub)blocks $A_{ij}$ may then have considerably different magnitudes---left untreated, this may again skew the above least squares formulation and lead to poor multigrid smoothing performance, e.g., because nearby elements are smoothed at different rates.
\end{itemize}
Fortunately, these aspects can be treated via a simple balancing process. In a form analogous to \cref{eq:blockform}, the $i$th diagonal block of $A$ is equal to 
\begin{equation} \label{eq:blockformdiag} \begin{pmatrix} \mathcal{A}_{ii} & \mathcal{G}_{ii} \\ \mathcal{D}_{ii} & 0 \end{pmatrix}. \end{equation}
Our balancing approach performs a simple diagonal pre- and post-scaling of $A$ in order to calibrate the norms of the (newly-scaled) blocks in \cref{eq:blockformdiag}. There are many possible approaches;\footnote{For example, one could simply compute the Frobenius norms of the matrices appearing in \cref{eq:blockformdiag} and use that to determine a diagonal scaling. This approach is easy, except for some edge cases concerning $\polydeg = 1$ for which it is possible that $\mathcal{D}_{ii} = 0$ and $\mathcal{G}_{ii} = 0$ on a small subset of mesh elements (depending on boundary conditions, numerical flux choices, etc.). We instead opted for an approach that works in all settings, for all $\polydeg$, and this led to the method described here based on the scaling of the broken gradient operator.} here we apply a scheme based on the block-diagonal broken gradient operator $\bg = (\bg_1, \ldots, \bg_d) : V_\polydeg \to V_\polydeg^d$ and how it factors into the operators of \cref{eq:blockformdiag}, as follows. We define the ``scale'' $\mathcal{S}_{\mathcal{A},i}$ of $\mathcal{A}_{ii}$ as
\[ \mathcal{S}_{\mathcal{A},i} := \Bigl(d\,\sum_{\ell = 1}^d \bigl\| (\bg_{\ell,ii})^\trans M_{\mu,ii}\,(\bg_{\ell,ii}) \bigr\|_F^2 \Bigr)^{\frac12}. \]
This definition reflects that, except for lifting operators, the viscous operator $\mathcal A$ is roughly $d$ copies of a $\mu$-weighted Laplacian operator. %
Similarly, we define the scale $\mathcal{S}_{\mathcal{D},i}$ of $\mathcal{D}_{ii}$ as
\[ \mathcal{S}_{\mathcal{D},i} := \Bigl(\sum_{\ell = 1}^d \bigl\| M_{ii}\,\bg_{\ell,ii} \bigr\|_F^2 \Bigr)^{\frac12}. \]
Using these scaling factors, let $\alpha \in \mathcal{V}_h$ denote a global scaling vector such that $\alpha_{i,\vu,\ell,k} := (S_{\mathcal{A},i})^{-1/2}$ and $\alpha_{i,p,k} := (S_{\mathcal{A},i})^{1/2} (S_{\mathcal{D},i})^{-1}$ for all $i,\ell,k$. The calibration of the blocks in \cref{eq:blockformdiag} is then implemented via the rescaled Stokes operator $\tilde A := \diag(\alpha) A \diag(\alpha)$; in particular, all of the newly-scaled subblocks have unit-magnitude Frobenius norm.

This balancing approach is applied as follows. Define $\tilde A := \diag(\alpha) A \diag(\alpha)$ and invoke the least squares methods discussed in \cref{sec:smoother1} on $\tilde A$ (instead of $A$) thereby computing the approximate inverse $\tilde Q$ satisfying $\tilde A \tilde Q \approx \mathbb{I}_\omega$. The true approximate inverse of $A$ is then defined by $Q := \diag(\alpha) \tilde Q \diag(\alpha)$ so that $A Q \approx \mathbb{I}_\omega$. One may view this process as a simple diagonal preconditioning of $A$ to treat high-contrast viscosity coefficients as well as the mesh-dependent scaling of the viscous stress, pressure gradient, and velocity divergence operators. One may also view this process as augmenting the weighted least squares problem via additional weighting factors that ``cancel out'' the scaling effects of $\mu$ and $h$.

\subsection{Unsteady Stokes problems}
\label{sec:unsteady}

Thus far we have focused on steady-state Stokes problems. For unsteady/nonstationary Stokes problems, the introduction of a density-weighted temporal derivative term can alter the characteristics of the Stokes system and therefore the multigrid solver. The governing equations consist of keeping the interfacial jump conditions \cref{eq:govern2} and boundary conditions \cref{eq:govern3} but amending \cref{eq:govern1} to include a density-weighted term that would arise from, e.g., a time-stepping method applied to the incompressible Navier-Stokes equations whereby the advection term is treated explicitly and the viscous term implicitly. The amended form reads 
\begin{equation} \label{eq:govern4} \left. \begin{aligned} \frac{\rho_i}{\dt} \vu - \nabla \cdot \bigl(\mu_i (\nabla \vu + \gamma\,\nabla \vu^\trans) \bigr) + \nabla p &= \vf \\ -\nabla \cdot \mathbf \vu &= f \end{aligned} \right\} \text{ in } \Omega_i, \end{equation}
where $\rho_i > 0$ is the density of phase $i$ and $\dt > 0$ is a parameter proportional to the time step of a temporal integration method. The corresponding LDG discretisation is a straightforward modification to the steady-state problem and leads to 
\begin{equation} \label{eq:blockformunsteady} \begin{pmatrix} \tfrac{1}{\dt} M_\rho + \mathcal{A} & \mathcal{G} \\ \mathcal{D} & 0 \end{pmatrix} \begin{pmatrix} \vu_h \\ p_h \end{pmatrix} = \begin{pmatrix} {\mathbf b}_{\vu} \\ b_p \end{pmatrix}\!, \end{equation}
where $M_\rho$ is a block diagonal $\rho$-weighted mass matrix. Beyond this simple modification, what makes the unsteady problem subtle is the relative strengths of the viscous operator and the newly-added density-weighted term \cite{flame,inferno}. When $\mu \dt/\rho$ is sufficiently large, we essentially have a steady-state Stokes system with a small $\rho/\dt$-weighted identity shift added to the viscous operator; in this case, a good steady-state multigrid Stokes solver can be effective. Conversely, when $\rho/(\mu \dt)$ is sufficiently large, \cref{eq:blockformunsteady} approximately reduces to a Poisson problem for pressure (written in flux-form), closely related to Chorin's projection method for incompressible fluid dynamics \cite{Chorin}. The role of the pressure variable $p$ changes in these two extremes, from being a Lagrange multiplier enforcing the divergence constraint (in the former), to being the primary solution variable in an elliptic interface problem (in the latter). Crucially, note also that the relative strengths of these operators can change across the levels of a multigrid hierarchy---on highly-refined meshes, viscous effects may dominate; on coarse grids, the density-weighted term may dominate.

These aspects can be treated via a refinement to the balancing process, together with a blending procedure applied to the smoother parameters $(\zetaI,\zetaII,\omegaI,\omegaII)$, as follows. 
\begin{itemize}
\item The ``scale'' of the temporal derivative term on the $i$th element is defined via 
\[ \mathcal{S}_{\rho,i} := \tfrac{1}{\dt} \| M_{\rho,ii} \|_F. \]
In the balancing process of \cref{sec:balance}, the scaling of the top-left block of \cref{eq:blockformunsteady} is defined as $(S_{\mathcal{A},i}^2 + \mathcal{S}_{\rho,i}^2)^{1/2}$, and we appropriately redefine the scaling vector $\alpha$.
\item The smoother parameters are blended via a {harmonic-weighted convex combination}. In one extreme, suppose that $\mu \neq 0$ and $\rho \equiv 0$ leads to a particular set of optimised smoother parameters $(\zetaI^0,\zetaII^0,\omegaI^0,\omegaII^0)$; in the other extreme, suppose that $\mu \equiv 0$ and $\rho \neq 0$ leads to optimised smoother parameters $(\zetaI^\infty,\zetaII^\infty,\omegaI^\infty,\omegaII^\infty)$. (The superscripts denote the effective Reynolds number.) Then, to build the smoother $Q_i$ on element $i$, for any combination of $\rho > 0$ and $\mu > 0$, we apply the following smoother parameters:
\[ (\zetaI,\zetaII,\omegaI,\omegaII) = (1 - \lambda_i) \bigl(\zetaI^0,\zetaII^0,\omegaI^0,\omegaII^0\bigr) + \lambda_i \bigl(\zetaI^\infty,\zetaII^\infty,\omegaI^\infty,\omegaII^\infty\bigr) \]
where
\[ \lambda_i := \frac{\mathcal{S}_{\rho,i}}{\mathcal{S}_{\rho,i} + \mathcal{S}_{\mathcal{A},i}}. \]
Here, $\lambda_i \in [0,1]$ is a blending amount that depends on the local, mesh-dependent relative strengths of the density-weighted temporal operator and the viscosity-weighted viscous operator. Note that the blending amount can vary across the multigrid hierarchy, and even between elements on the same level.
\end{itemize}

\subsection{Smoother construction and application}

\begin{figure}[!t]%
\centering%
\begin{minipage}{0.65\textwidth}%
\begin{algorithm}[H]
	\caption{\sffamily\small Construction of $Q$, given smoother parameters $(\zetaI^0,\zetaII^0,\omegaI^0,\omegaII^0)$ and $(\zetaI^\infty,\zetaII^\infty,\omegaI^\infty,\omegaII^\infty)$ for the extremal cases of $\{\mu \neq 0, \rho \equiv 0\}$ and $\{\mu \equiv 0, \rho \neq 0\}$, respectively.}
	\begin{algorithmic}[1]
            \For{every element $i$ (possibly in parallel)}
                \State $\mathcal{S}_{\mathcal{A},i} := \Bigl(d\,\sum_{\ell = 1}^d \bigl\| (\bg_{\ell,ii})^\trans M_{\mu,ii}\,(\bg_{\ell,ii}) \bigr\|_F^2 \Bigr)^{\frac12}.$
                \State $\mathcal{S}_{\mathcal{D},i} := \Bigl(\sum_{\ell = 1}^d \bigl\| M_{ii}\,\bg_{\ell,ii} \bigr\|_F^2 \Bigr)^{\frac12}.$
                \State $\mathcal{S}_{\rho,i} := \begin{cases} 0 & \text{for steady-state Stokes problems,} \\ \tfrac{1}{\dt} \| M_{\rho,ii} \|_F & \text{for unsteady Stokes problems}. \\  \end{cases}$
            \EndFor
            \State Define the rescaling vector $\alpha \in \mathcal{V}_h$ such that $\forall i,\ell,k$,
            \[ \begin{aligned} \alpha_{i,\vu,\ell,k} &:= (S_{\mathcal{A},i}^2 + \mathcal{S}_{\rho,i}^2)^{-1/4}, \\ \alpha_{i,p,k} &:= (S_{\mathcal{A},i}^2 + \mathcal{S}_{\rho,i}^2)^{1/4} (S_{\mathcal{D},i})^{-1}. \end{aligned} \]
            \State Define the balanced Stokes operator $\tilde A := \diag(\alpha) A \diag(\alpha)$.
            \For{every element $i$ (possibly in parallel)}
                \State $\lambda_i := \mathcal{S}_{\rho,i} / (\mathcal{S}_{\rho,i} + \mathcal{S}_{\mathcal{A},i})$.
                \State $(\zetaI,\zetaII,\omegaI,\omegaII) := (1 - \lambda_i) (\zetaI^0,\zetaII^0,\omegaI^0,\omegaII^0) + \lambda_i (\zetaI^\infty,\zetaII^\infty,\omegaI^\infty,\omegaII^\infty).$
                \State Define $\zeta_i \in \mathcal{V}_h$ such that $\zeta_i \equiv 0$, except for
                \[ \begin{aligned} \zeta_{i;i,\vu,\ell,k} &:= \zetaI && \forall \ell,k, \\ \zeta_{i;i,p,k} &:= \zetaII &&\forall k, \\ \zeta_{i;j,\vu,\ell,0} &:= 1 &&\forall \ell, j \neq i. \end{aligned} \]
                \State Define $\omega_i \in \mathcal{V}_h$ such that $\omega_i \equiv 0$, except for
                \[ \begin{aligned} \omega_{i;i,\vu,\ell,k} &:= \omegaI &&\forall \ell,k, \\ \omega_{i;i,p,k} &:= \omegaII &&\forall k. \end{aligned} \]
                \State Solve the least squares minimisation problem \label{alg:Qi}
                \[ \tilde{Q}_i := \argmin_{\mathsf{Q} \in \mathbb{R}^{n \times n}}~\bigl\| \diag(\zeta_i) \bigl( \tilde{A}_{:,i}\,\mathsf{Q} - \diag(\omega_i)\, \mathbb{I}_{:,i} \bigr) \bigr\|_F^2. \]
            \EndFor
		\State{\textbf{return} $Q := \diag(\alpha) \tilde{Q} \diag(\alpha)$.}
	\end{algorithmic}%
	\label{algo:Q}%
\end{algorithm}%
\end{minipage}%
\end{figure}

Assembling the above steps, \Cref{algo:Q} details the overall approach to construct the block-diagonal smoother $Q$. We note several aspects:
\begin{itemize}
\item The implementation of \cref{algo:Q} is perhaps simpler than it may first appear. The main step is that of assembling and solving the least squares problems. In essence, assembly involves visiting every nonzero block in the $i$th block column of $\tilde A$---upon visiting the diagonal block, rescale the rows arising from the stress divergence (resp., velocity constraint) operator by $\zetaI$ (resp., $\zetaII$) and append them to a temporary matrix; upon visiting an off-diagonal block, keep only the rows arising from the stress divergence operator associated with constant-mode lowest-degree modal coefficients, and append them to the same temporary matrix. The resulting matrix is of size $m \times n$ where $n = d(\polydeg + 1)^d + \polydeg^d$ and $m = n + \ell d$, where $\ell$ is the number of off-diagonal blocks (often just four in 2D or six in 3D). A similar temporary matrix is built to define the right-hand side of the least squares problem, i.e., an $m \times n$ matrix with either $\zetaI \omegaI$ or $\zetaII \omegaII$ appearing on the diagonal entries. The Frobenius norm minimisation problem (line \ref{alg:Qi}) translates to a standard least squares problem with $n$ right-hand side vectors; in this work we utilise a standard QR method from \textsc{lapack} for its solution.
\item Each individual least squares problem involves a nearly square matrix with computational cost similar to inverting a diagonal block of $A$; further comments on the size of these matrices are given in the concluding remarks. Moreover, each least squares problem is well-conditioned and almost always full-rank.\footnote{The exception is when a single mesh element spans the entire domain, plus some idiosyncratic cases involving $\polydeg = 1$; in these cases the rank deficiency is straightforward to handle, simply by adopting the minimum-norm minimum-residual least squares solution.}
\item Note that every diagonal block of $Q$ can be (pre)computed independently of (and possibly concurrently with) every other diagonal block; i.e., the algorithm is trivially parallelisable.
\end{itemize}

We now discuss the application of $Q$ as a multigrid smoother, analogous to multicoloured block Gauss-Seidel. Let the input approximation be denoted by $\xin \in \mathcal{V}_h$ and the output by $\xout \in \mathcal{V}_h$. In general, $Q$ is asymmetric, so there are two alternatives to the Gauss-Seidel sweep of \cref{eq:gssweep}: the first uses $Q$,
\begin{equation} \label{eq:Qpre} \xout_i := \xin_i + Q_i \bigl(b_i - \textstyle \sum_{j} A_{ij} \mathring{x}_{j,i} \bigr) \quad \text{where} \quad \mathring{x}_{j,i} := \begin{cases} \xout_j & \text{if } j < i, \\ \xin_j & \text{otherwise,} \end{cases} \end{equation}
while the second uses the transpose of $Q$,
\begin{equation} \label{eq:Qpost} \xout_i := \xin_i + Q_i^\trans \bigl(b_i - \textstyle \sum_{j} A_{ij} \mathring{x}_{j,i} \bigr) \quad \text{where} \quad \mathring{x}_{j,i} := \begin{cases} \xout_j & \text{if } j < i, \\ \xin_j & \text{otherwise,} \end{cases} \end{equation}
where as before, ``$j < i$'' means block $j$ is updated before block $i$. One could therefore apply any combination of \cref{eq:Qpre} or \cref{eq:Qpost} as a pre-smoother or post-smoother. Numerical experiments indicate that the fastest multigrid solvers arise from applying \cref{eq:Qpre} as a pre-smoother and \cref{eq:Qpost} as a post-smoother---consequently, the same choice is used throughout the numerical experiments presented below; see also the results of \cite{inferno} that arrive at a similar conclusion for smoothers based on general SAI methods. Regarding the element sweep order, we consider here a multicoloured approach. Specifically, in a setup phase, a graph-colouring algorithm is applied to the connectivity graph associated with the block sparsity pattern of the Stokes operator $A$, thereby assigning a colour to each element.\footnote{For example, one could apply the \textsc{DSatur} graph-colouring algorithm \cite{lewis2015guide}.} (On Cartesian grids with one-sided numerical fluxes, this recovers the well-known red-black colouring associated with the 5-point (2D) or 7-point (3D) Laplacian stencil.) In the pre-smoothing sweep of \cref{eq:Qpre} and post-smoothing sweep of \cref{eq:Qpost}, every element of the same colour is processed (possibly in parallel), before moving onto the next colour.\footnote{In particular, the pre- and post-smoothing sweeps use the same colour ordering, a choice motivated by the findings of \cite{inferno}, where this configuration typically yielded faster convergence than colour reversal during post-smoothing. %
}

\subsection{Optimal smoothing parameters}
\label{sec:lfa}

To summarise the discussion so far, we have built a smoother analogous to multicoloured block Gauss-Seidel, taking into account key aspects of multiphase Stokes problems such as mesh-dependent operator scaling, as well as the application of the smoother in a multigrid V-cycle. The smoother has four tunable parameters: $\zetaI$ and $\zetaII$, controlling the relative strengths of velocity and pressure variables in the approximate inverses, and $\omegaI$ and $\omegaII$, controlling the damping/overrelaxation of velocity and pressure degrees of freedom. To conclude this section, we discuss the optimisation of these parameters.

Numerical experiments universally indicate that the fastest multigrid solvers are obtained when $\zetaI$ is sufficiently large ($\zetaI \gtrapprox 10$), beyond which the precise value of $\zetaI$ plays no further role in achieving the best performance; in other words, one can essentially make $\zetaI$ arbitrarily large. Since $\zetaI$ implements a weighting of the stress divergence operator on the diagonal block, when $\zetaI$ is sufficiently large the corresponding equations are nearly exactly satisfied by the least squares solution---intuitively, the smoothing effects of the Laplace-like viscous operator $\mathcal A$ should be prioritised. Note that this implies that the corresponding unknowns could be eliminated (essentially via the Schur complement of $\mathcal{A}_{ii}$), thereby reducing computational cost of the least squares solver, albeit with more algorithmic bookkeeping. On the other hand, in our numerical experiments we observed that the (pre)computation of the smoother $Q$ has negligible cost compared to the overall multigrid solver. Consequently, and mainly for algorithmic simplicity, we have kept the non-eliminated form and apply a large value of $\zetaI$. In particular, we universally set $\zetaI = 128$, across every parameter combination and every test problem, and exclude $\zetaI$ from further consideration.\footnote{The value of 128 carries no particular meaning, other than being one order magnitude greater than the aforementioned cutoff and not so excessive that it impacts numerical conditioning.}

To optimise for the remaining parameters $(\zetaII,\omegaI,\omegaII)$, here we adopt an approach based on two-grid local Fourier analysis \cite{mgbook1,mgbook2,doi:10.1137/19M1308669}, detailed further in \cref{app:lfa}. In essence, this procedure evaluates the performance of a two-level multigrid method, as applied to an infinite uniform Cartesian grid spanning all of $\R^d$, on an input function space spanned by functions of the form $\bm x \mapsto \exp(\iota \bm\theta \cdot \bm x/h)$, where $\iota^2 = -1$. %
The output is a ``two-grid LFA convergence factor'' $\rho(\zetaII,\omegaI,\omegaII)$, a predicted multigrid convergence rate such that the error in solving $Ax = b$ is reduced by a factor of $\rho$ per iteration. A corresponding search procedure then seeks out the set of parameters for which convergence is fastest, i.e., $\rho(\zetaII,\omegaI,\omegaII)$ is smallest. Given the nature of these parameters---specifically that $\omegaI$ and $\omegaII$ are damping/overrelaxation parameters, while the balancing algorithms are designed so that $\zetaII$ is unit magnitude---we expect a (local) minimum near the point $(1,1,1)$; the search algorithm applies a heuristic coordinate descent method to approximately find this minimum, denoted in the following by $\rho_\star$, and the corresponding tuned parameters $(\zetaII,\omegaI,\omegaII)$. {Note that two-grid LFA serves primarily as a guiding heuristic, since it optimises performance under the idealised assumptions of an infinite, constant-coefficient uniform grid and an exact coarse-grid solve. Nevertheless, the numerical experiments in \cref{sec:experiment} demonstrate that these LFA-tuned parameters transfer effectively to more complex settings, including nontrivial boundary conditions, high-contrast variable viscosity coefficients, and nonuniform meshes.}

\newcommand{\panel}[7]{\begin{minipage}{1.75in}\centering%
\small$\wp = #3$, \IfEqCase{#4}{{2}{$\Rey = \infty$}}[$\gamma = #4$]\\[0.5ex]%
\footnotesize$\mathsf{\ball = (#5)}$\\[0.25ex]%
\scriptsize$\mathsf{\cube = (#6)\!\nearrow\!(#7)}$\\%
\includegraphics[width=\textwidth]{#1}%
\end{minipage}}

\begin{figure}[!t]%
\centering%
\panel{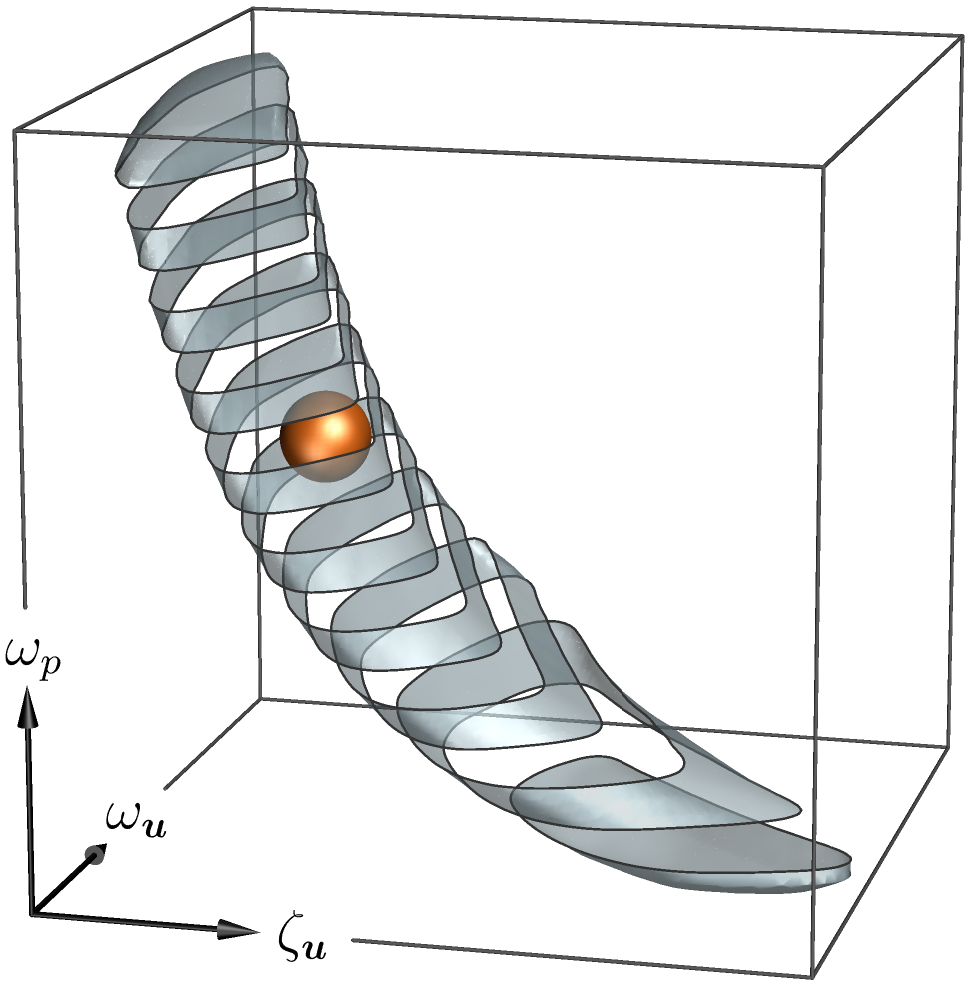}{2}{1}{0}{0.406,1.19,1.06}{0.23,1,0.24}{0.91,1.4,1.8}\quad%
\panel{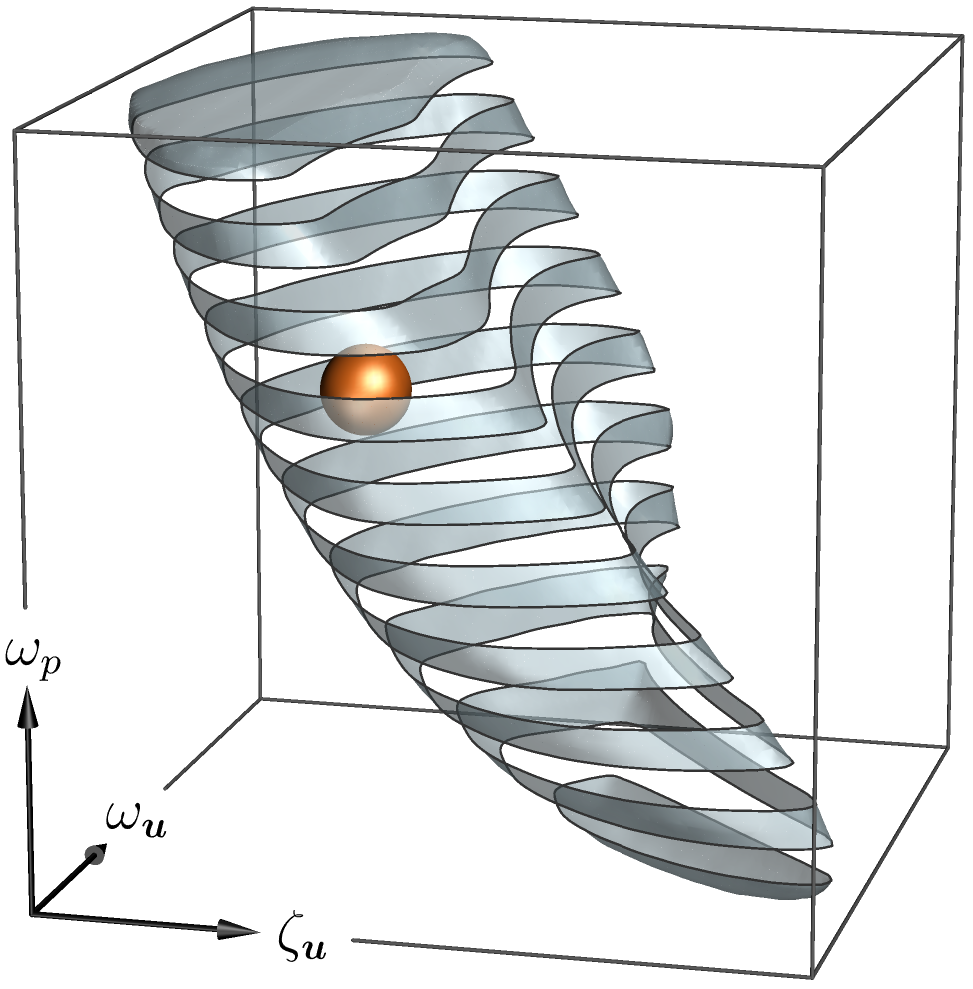}{2}{1}{1}{0.594,1.22,1.25}{0.44,1,0.53}{0.95,1.4,1.8}\quad%
\panel{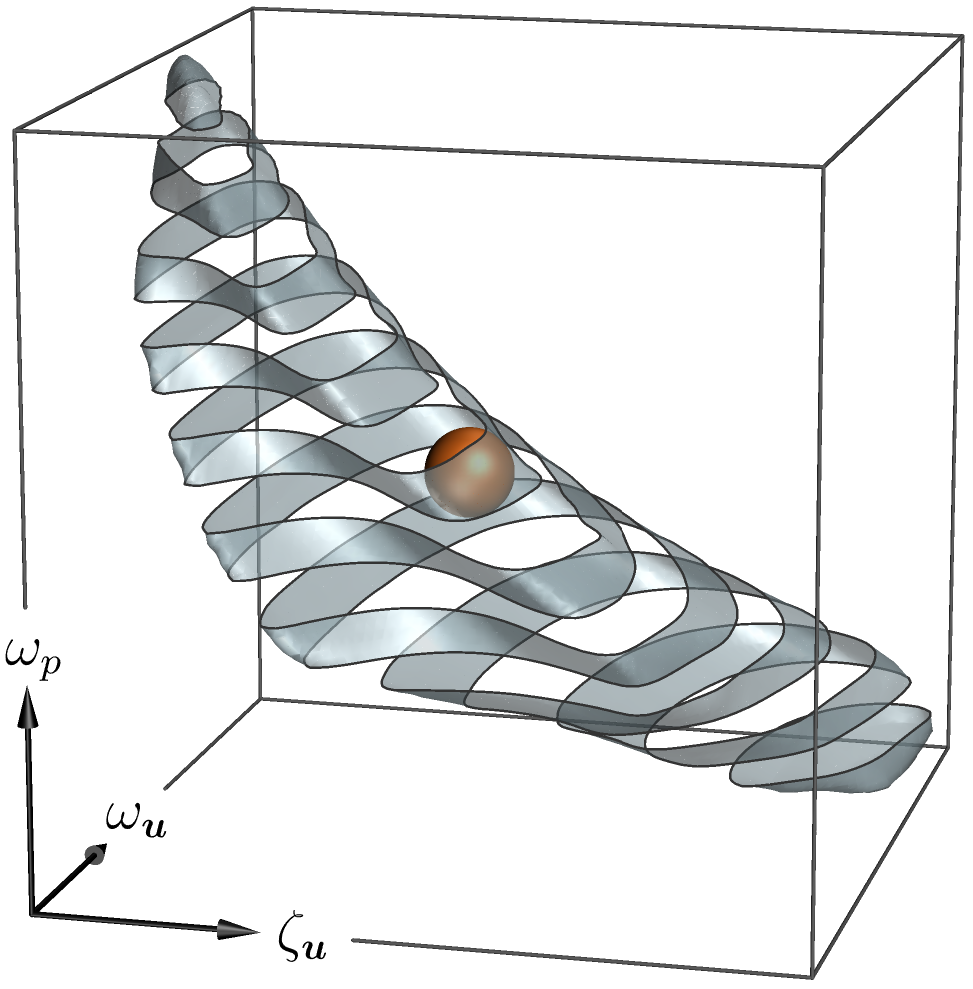}{2}{1}{2}{0.414,1,1.09}{0.31,0.95,0.74}{0.56,1,1.5}\\[1em]%
\panel{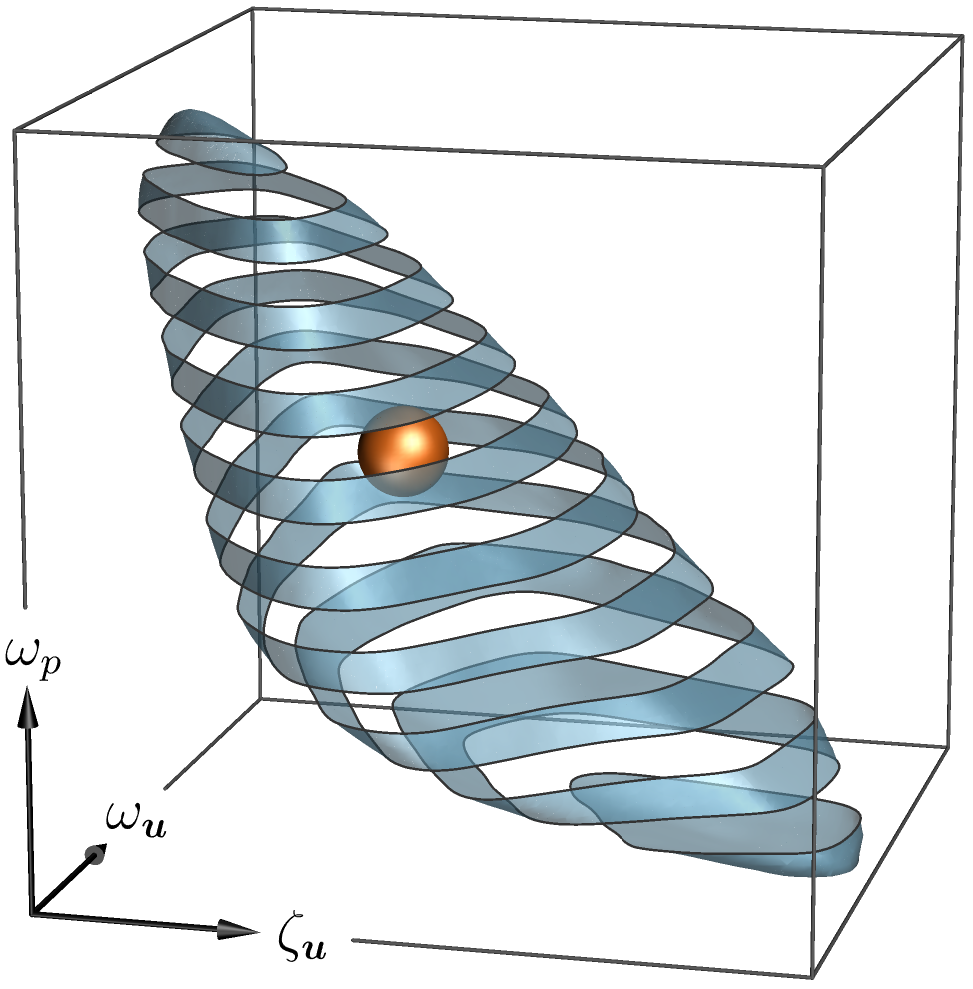}{2}{2}{0}{0.563,0.984,0.922}{0.39,0.88,0.73}{0.88,1.1,1.1}\quad%
\panel{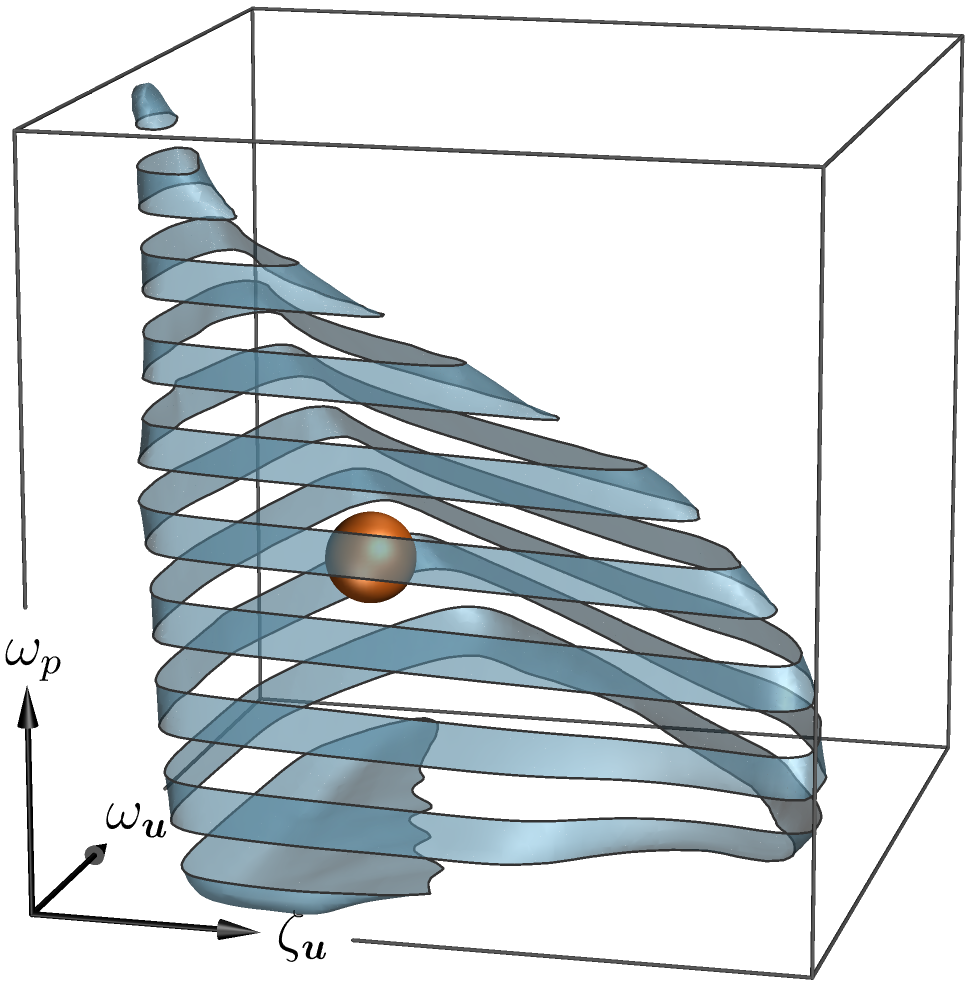}{2}{2}{1}{1.38,0.984,0.75}{0.49,0.93,0.62}{3.3,1.1,0.97}\quad%
\panel{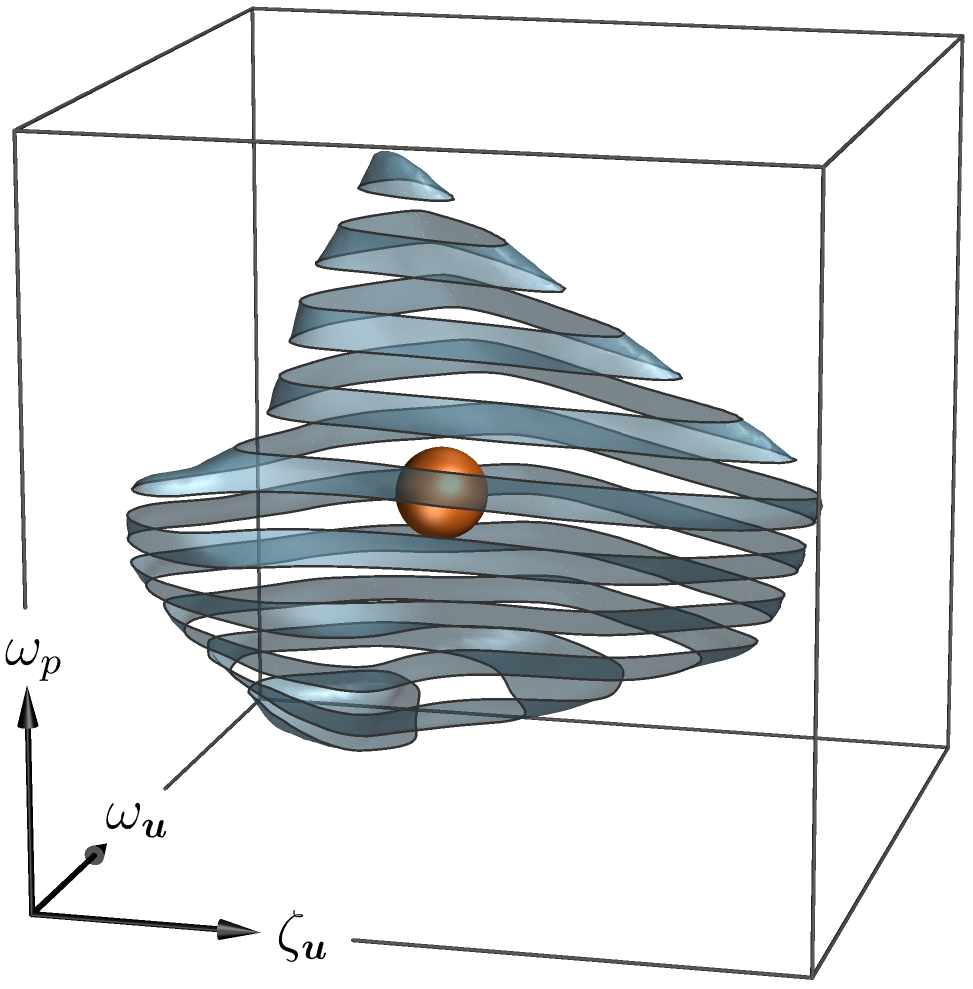}{2}{2}{2}{3.13,0.734,0.625}{2.1,0.65,0.48}{4.6,0.86,0.78}\\[1em]%
\panel{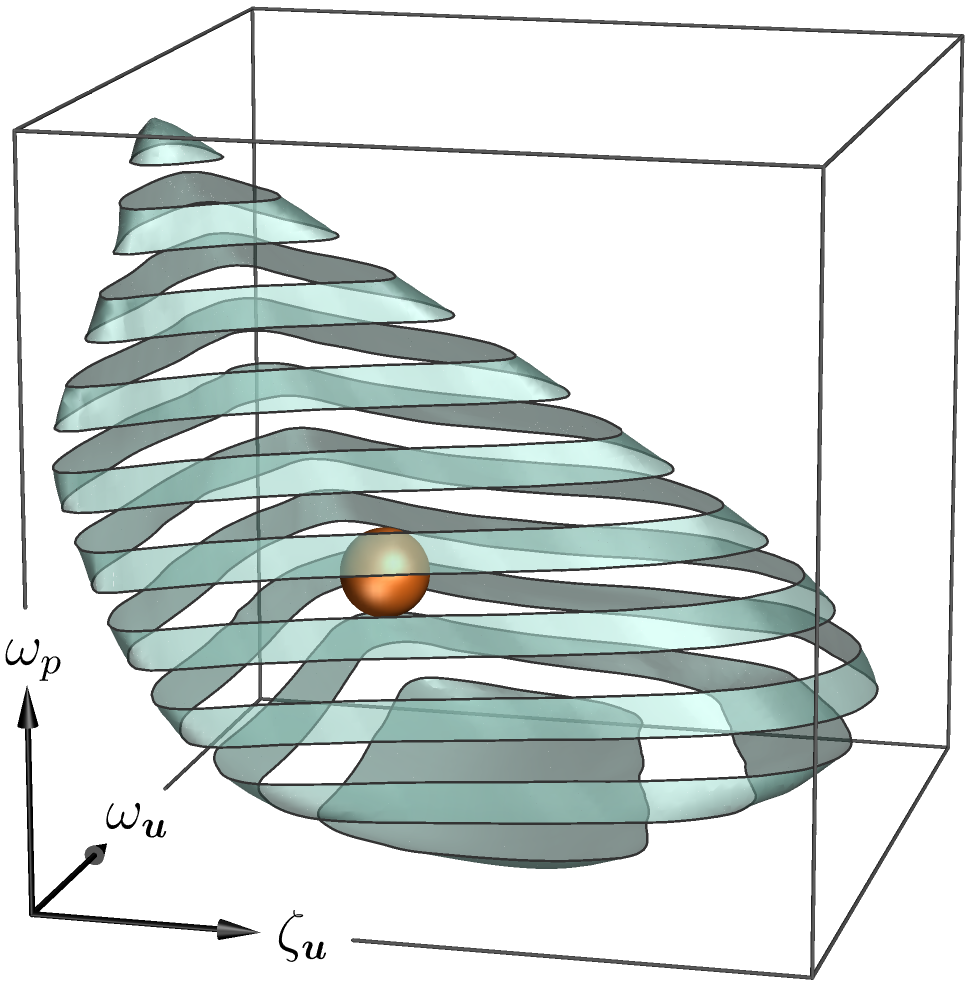}{2}{3}{0}{0.563,1.02,0.742}{0.32,0.94,0.66}{1.1,1.1,0.91}\quad%
\panel{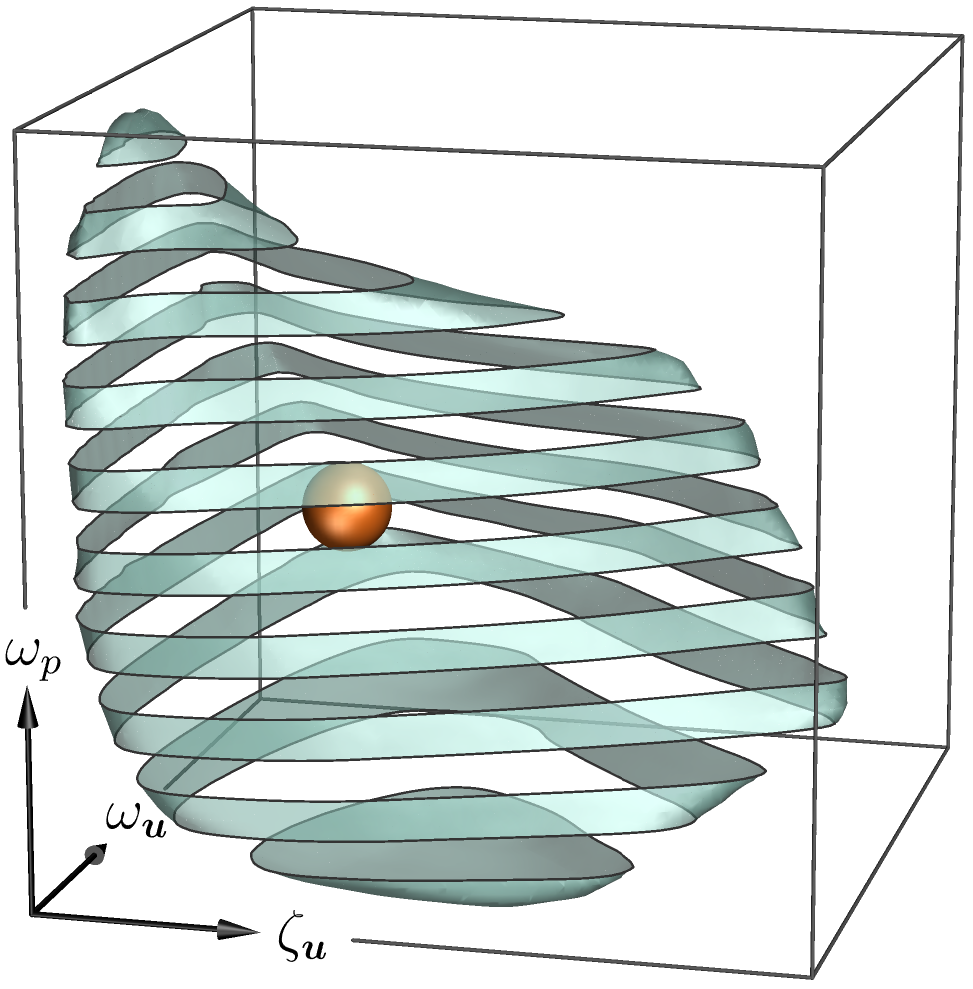}{2}{3}{1}{0.813,1.1,0.766}{0.35,1,0.63}{2.2,1.2,0.96}\quad%
\panel{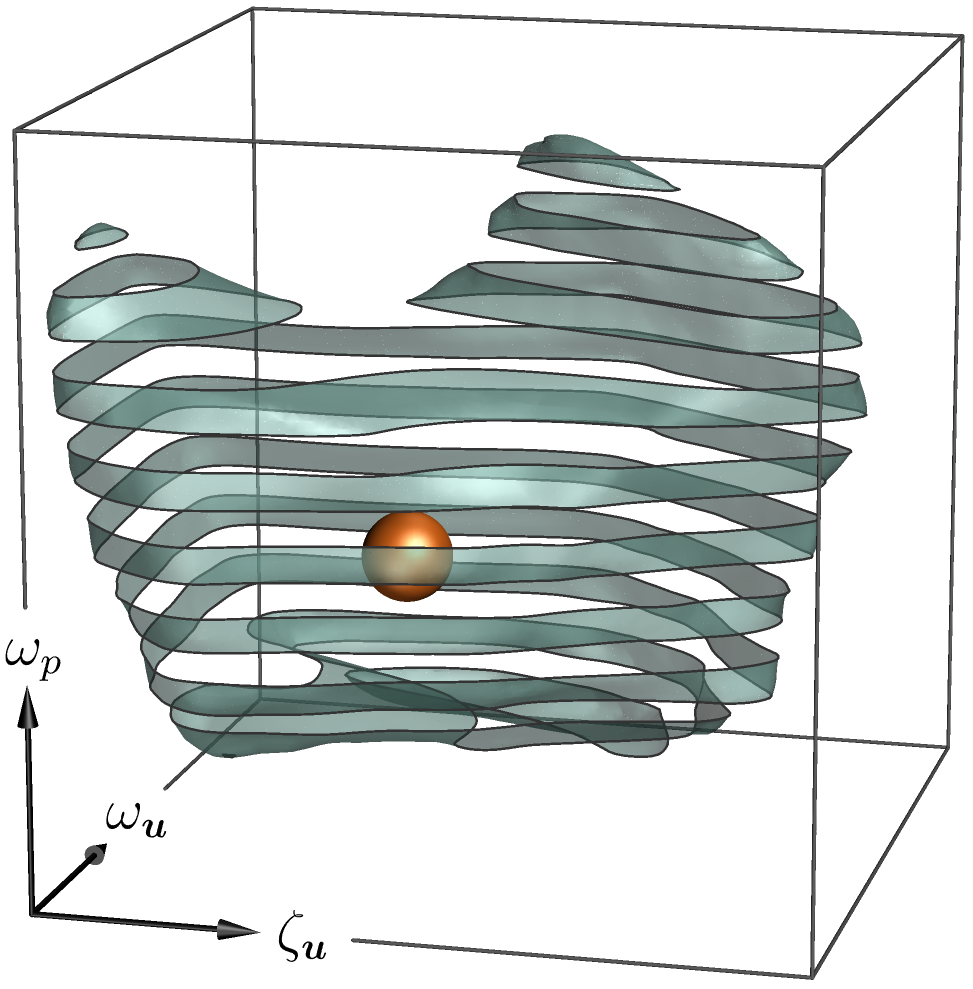}{2}{3}{2}{3.5,0.906,0.531}{2.7,0.79,0.42}{4.9,1,0.71}
\caption{Regions of near-optimal smoother parameters for 2D steady-state standard-form Stokes problems (left column), steady-state stress-form Stokes problems (middle column), and unsteady vanishing-viscosity Stokes problems (right column), for each $\polydeg \in \{1, 2, 3\}$. Here, ``near-optimal'' means the multigrid iteration count is at most 10\% above optimal, according to the predictions of two-grid local Fourier analysis, as outlined in \cref{sec:lfa} and detailed further in \cref{app:lfa}. The approximate centroid $\ball$ and bounding box $\cube$ of each region is indicated at the top of each panel; {see \cref{sec:lfa} for their definition.}}%
\label{fig:optimal2}%
\end{figure}

\begin{figure}[!t]%
\centering%
\panel{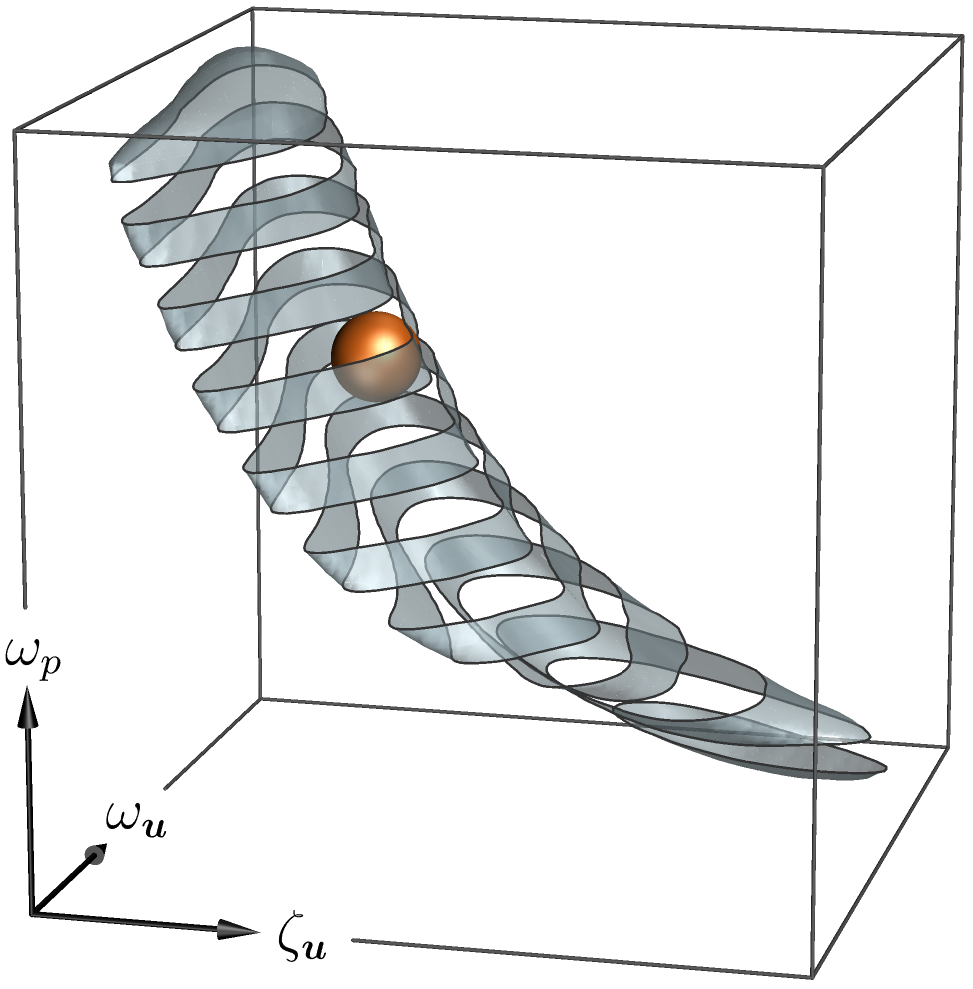}{3}{1}{0}{0.328,1.13,1.25}{0.21,0.81,0.34}{0.61,1.4,1.8}\quad%
\panel{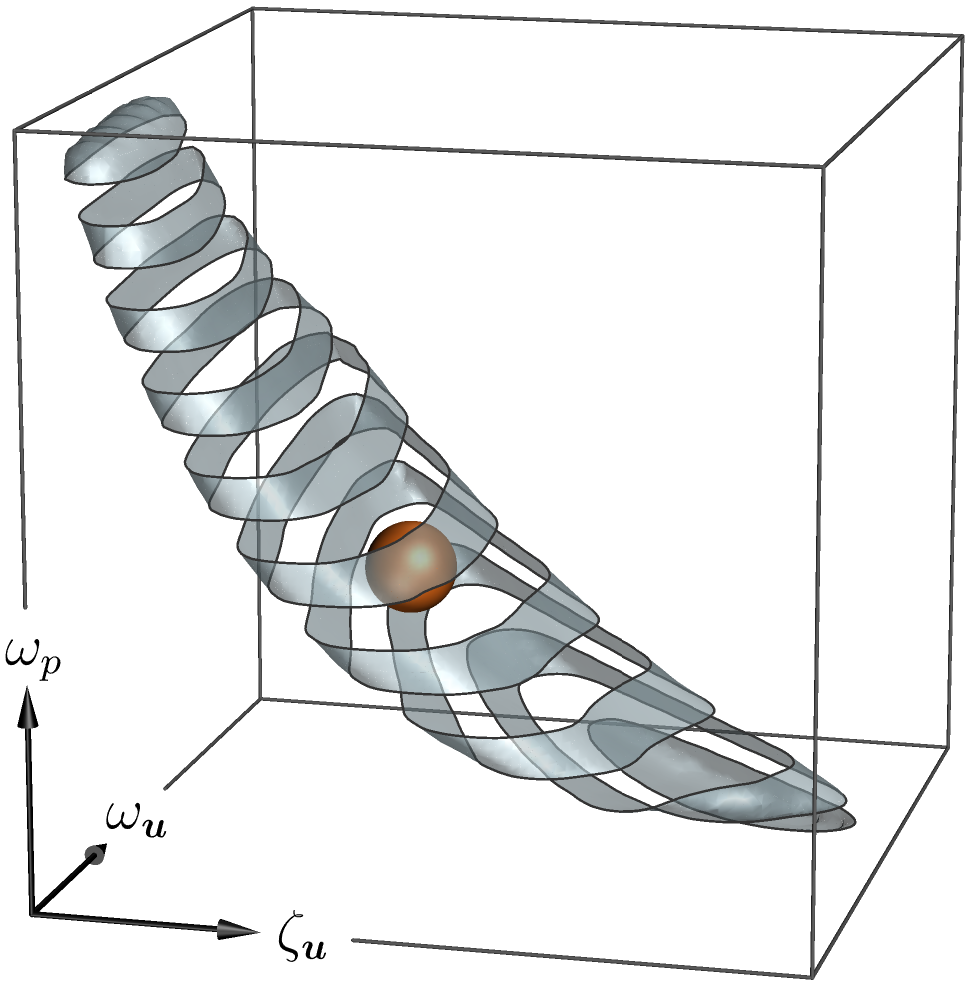}{3}{1}{1}{0.563,1.09,0.938}{0.36,0.92,0.43}{0.9,1.3,1.8}\quad%
\panel{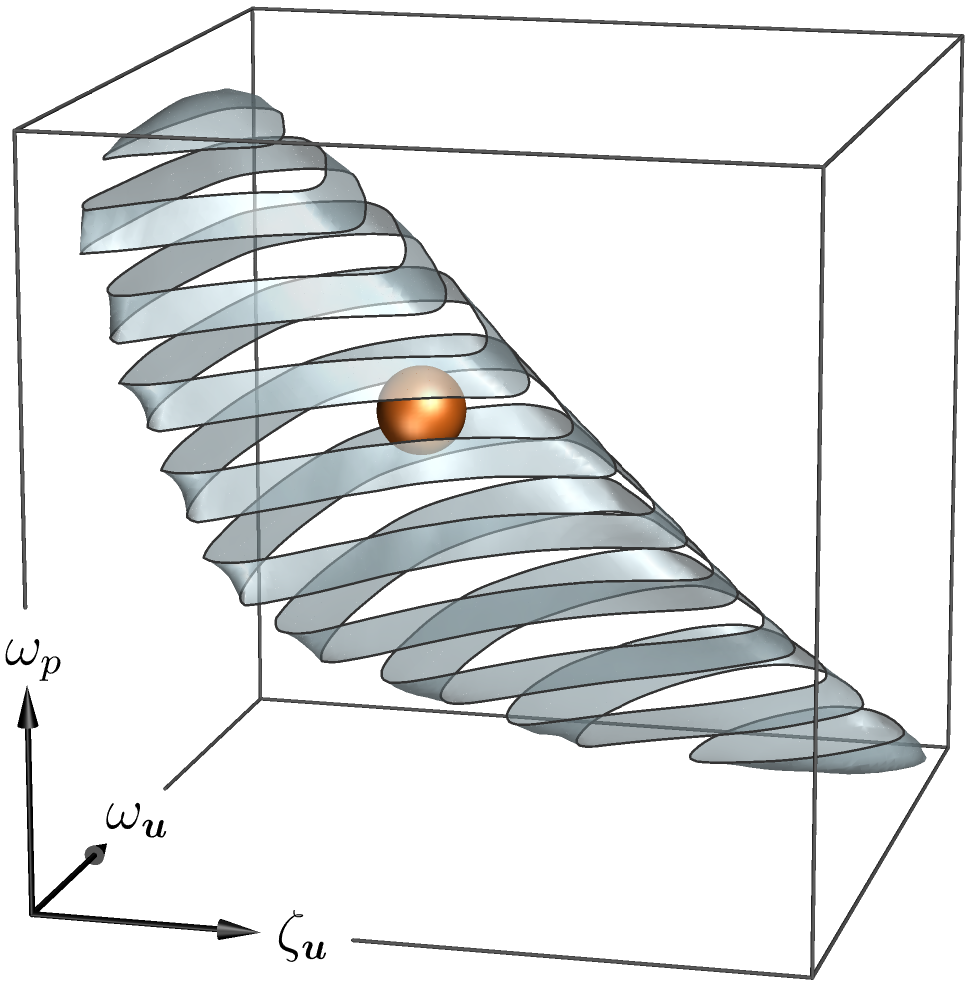}{3}{1}{2}{0.367,1.01,1.03}{0.32,0.95,0.81}{0.45,1,1.2}\\[1em]%
\panel{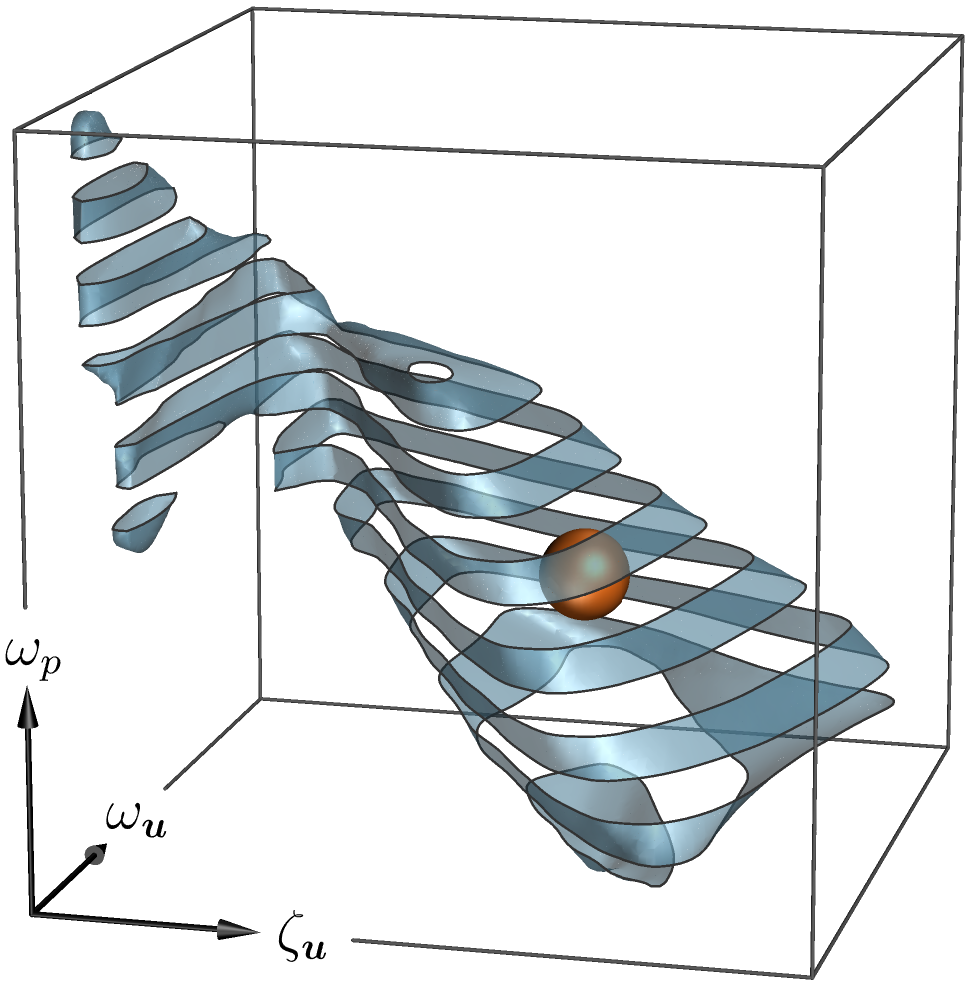}{3}{2}{0}{0.875,0.992,0.578}{0.28,0.9,0.43}{1.3,1.1,0.84}\quad%
\panel{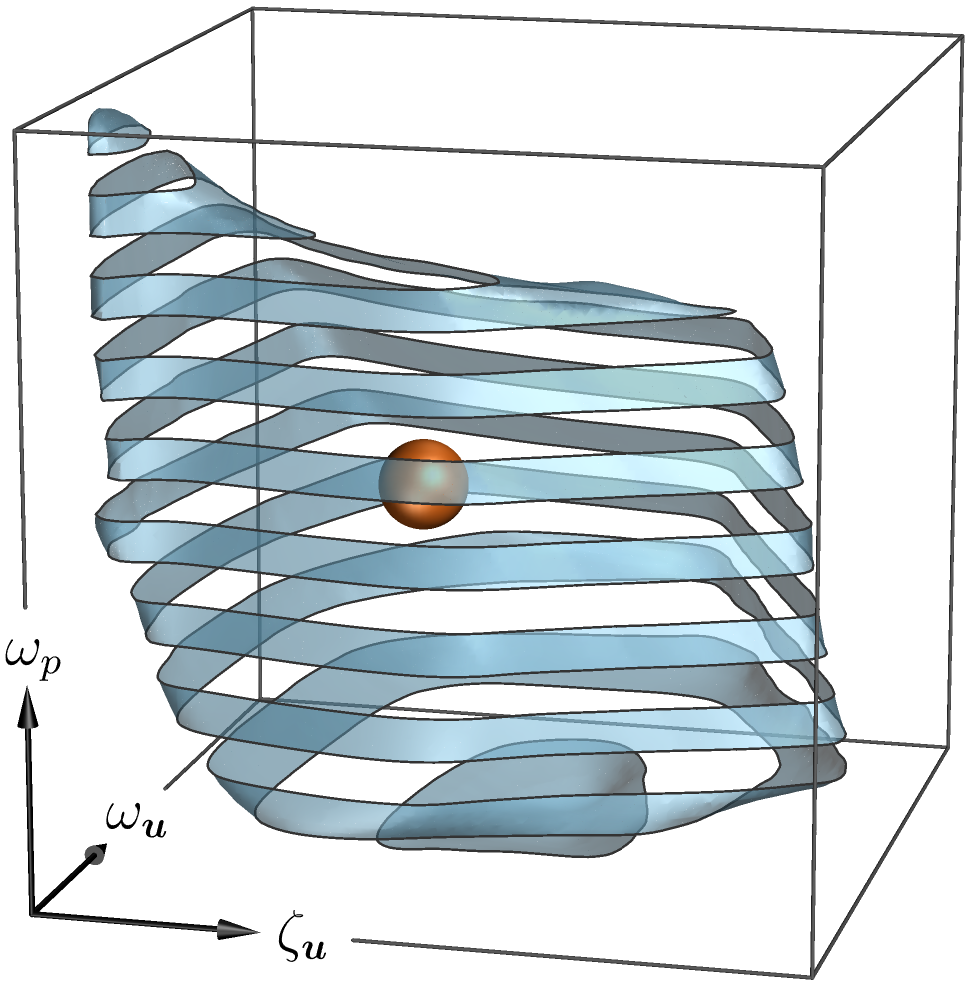}{3}{2}{1}{1,1.04,0.625}{0.32,0.96,0.43}{2.2,1.1,0.87}\quad%
\panel{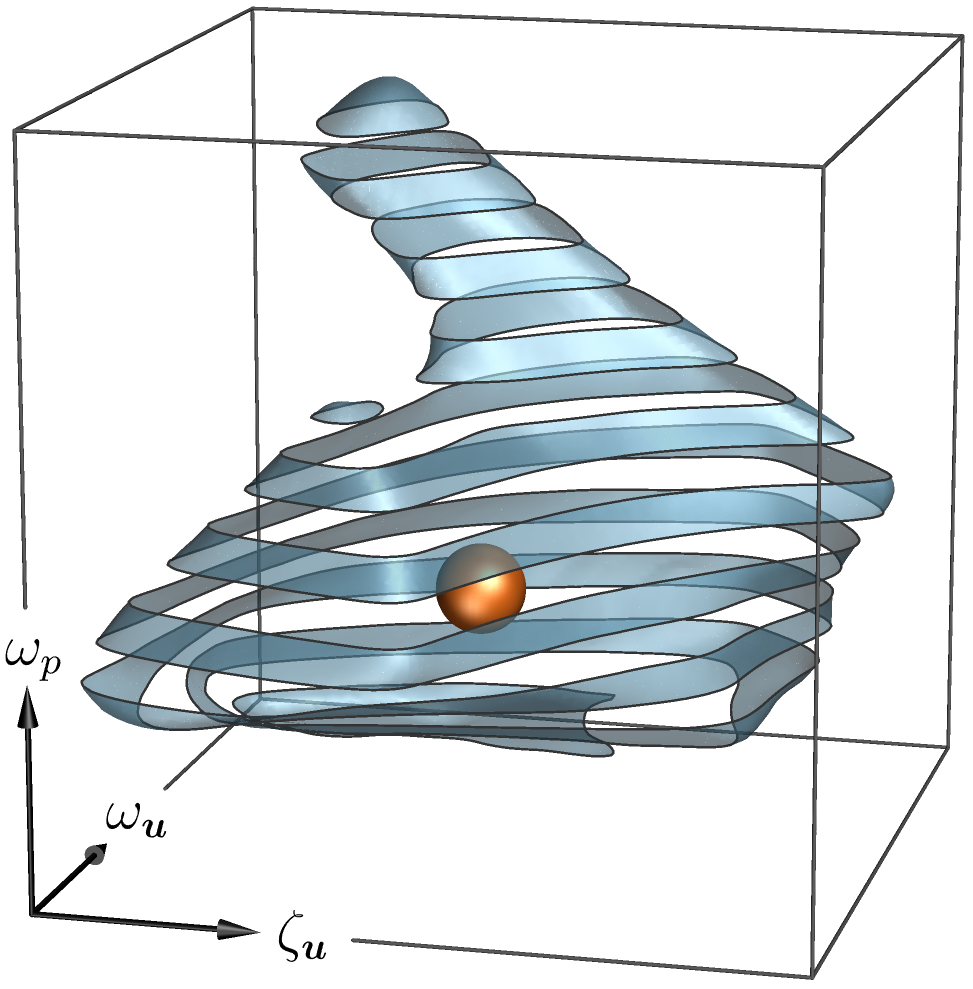}{3}{2}{2}{2.38,0.75,0.516}{1.9,0.63,0.44}{2.9,0.84,0.69}\\[1em]%
\panel{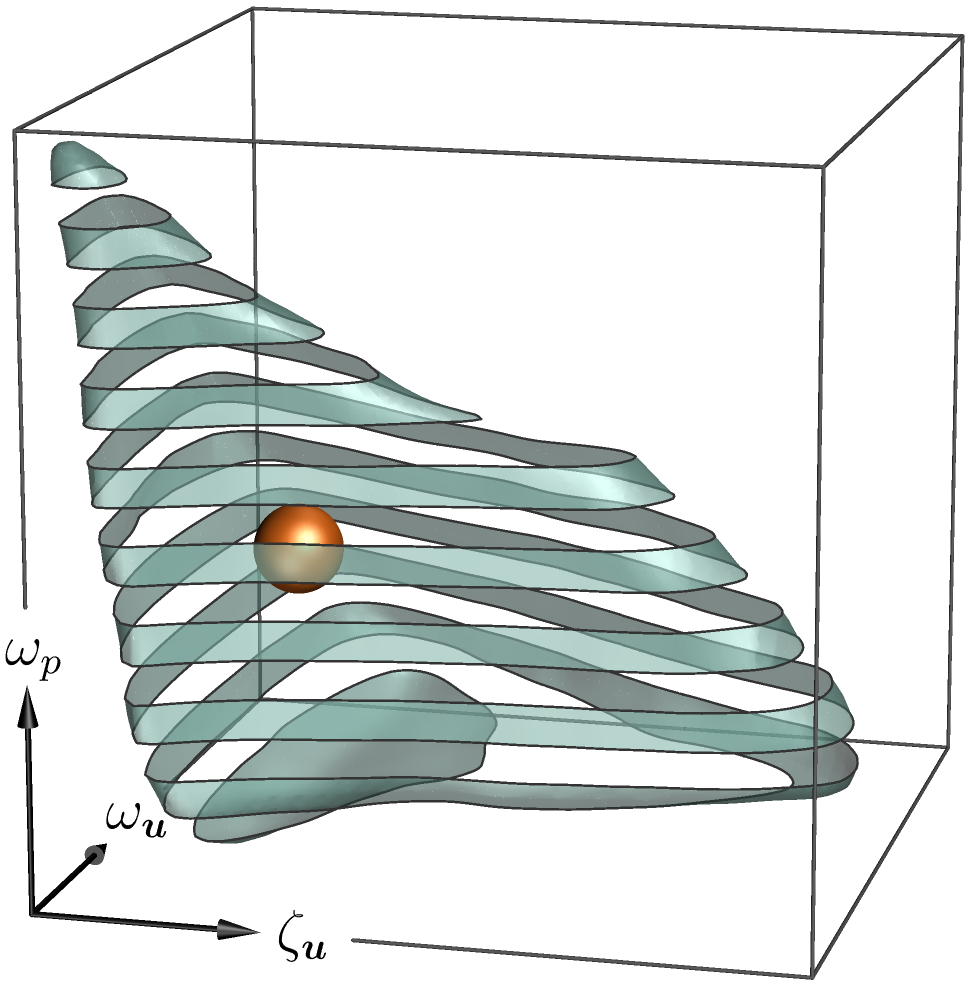}{3}{3}{0}{0.5,0.984,0.531}{0.25,0.93,0.43}{1.5,1,0.71}\quad%
\panel{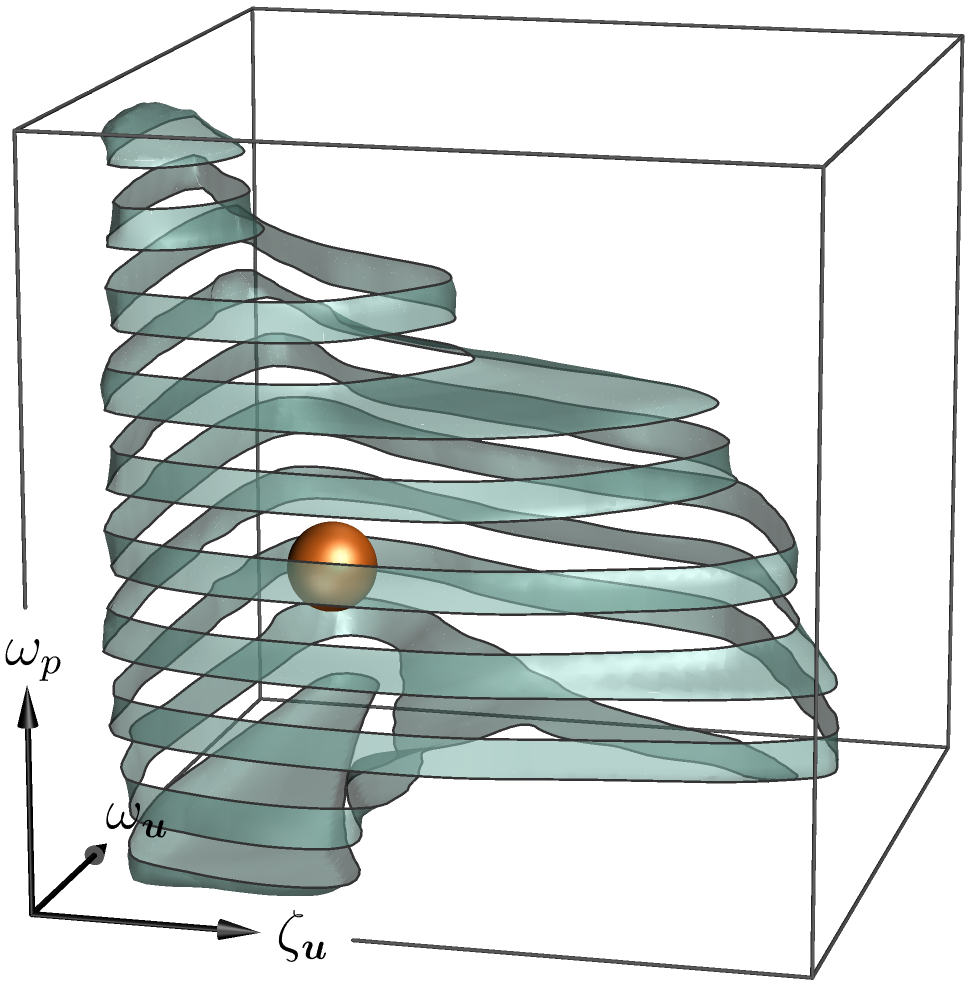}{3}{3}{1}{0.875,1.04,0.578}{0.37,0.98,0.51}{2.5,1.1,0.71}\quad%
\panel{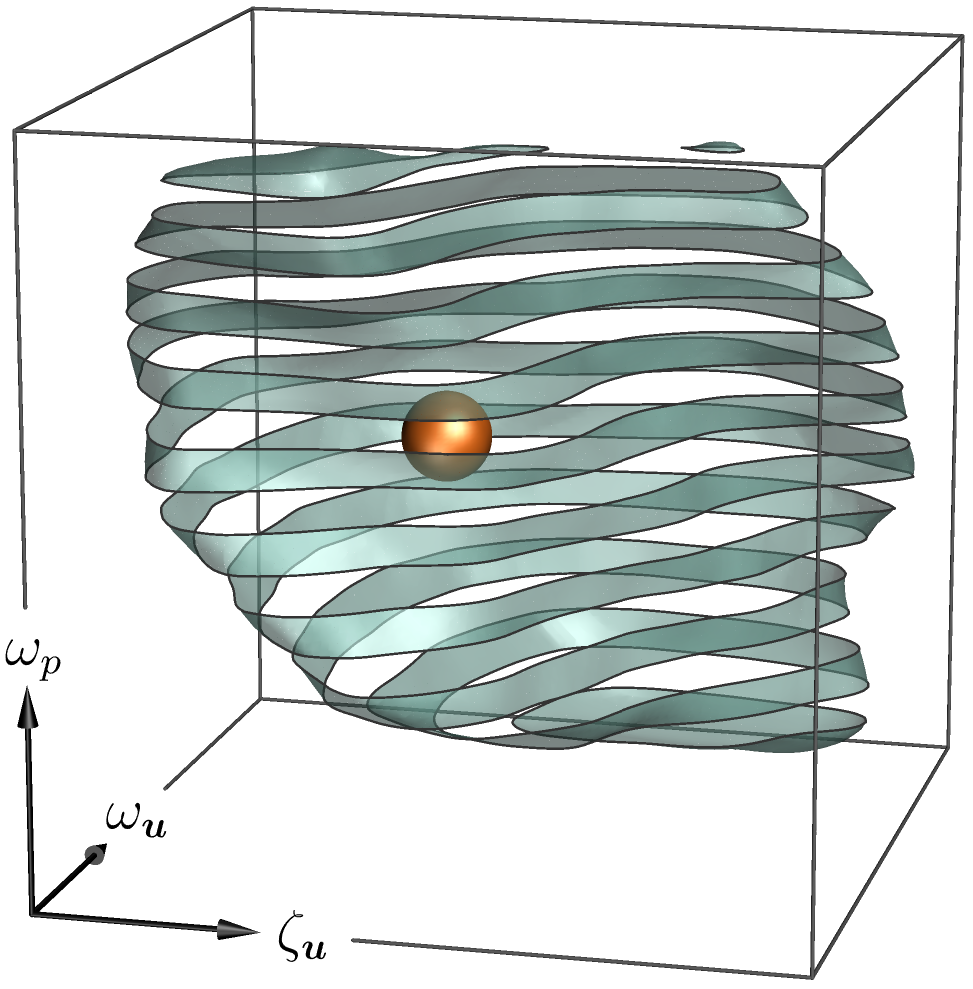}{3}{3}{2}{2.25,0.797,0.578}{1.7,0.71,0.45}{3.1,0.87,0.7}
\caption{Regions of near-optimal smoother parameters for 3D steady-state standard-form Stokes problems (left column), steady-state stress-form Stokes problems (middle column), and unsteady vanishing-viscosity Stokes problems (right column), for each $\polydeg \in \{1, 2, 3\}$. Here, ``near-optimal'' means the multigrid iteration count is at most 10\% above optimal, according to the predictions of two-grid local Fourier analysis, as outlined in \cref{sec:lfa} and detailed further in \cref{app:lfa}. The approximate centroid $\ball$ and bounding box $\cube$ of each region is indicated at the top of each panel; {see \cref{sec:lfa} for their definition.}}%
\label{fig:optimal3}%
\end{figure}

For every combination of polynomial degree $\polydeg \in \{1,2,3\}$ and spatial dimension $d \in \{2,3\}$, we consider three cases: (i) a steady-state standard-form Stokes problem, in which $\rho \equiv 0$, $\mu \equiv 1$, and $\gamma = 0$; (ii) a steady-state stress-form Stokes problem, in which $\rho \equiv 0$, $\mu \equiv 1$, and $\gamma = 1$; and (iii) an unsteady Stokes problems of infinite Reynolds number, in which $\rho \equiv 1$ and $\mu \equiv 0$ ($\gamma$ thus irrelevant). For each of these 18 configurations, we apply the two-grid LFA method and at the same time, evaluate $\rho$ in a neighbourhood of the minimum $\rho_\star$. The results of this study are shown in \cref{fig:optimal2,fig:optimal3} for the 2D and 3D cases, respectively:
\begin{itemize}
    \item The illustrated surfaces correspond to the level set $\{\rho(\zetaII,\omegaI,\omegaII) = \rho_\star^{1/1.1}\}$ and indicate a region of near-optimal parameters. Any set of parameters inside this region yields a two-grid LFA-predicted multigrid convergence rate no more than $\smash{\rho_\star^{1/1.1}}$; equivalently, they yield a multigrid iteration count at most 10\% above the LFA-predicted optimal.%
    \item The bounding box of the near-optimal region provides an indication of how finely-tuned the parameters need to be. For example, the parameter $\zetaII$ often has a relatively wide window, whereas $\omegaI$ is often more tightly controlled around unity. This bounding box is indicated at the top of each figure; e.g., ``{\small$\cube = (0.4,0.9,0.7)\!\nearrow\!(0.9,1.1,1.1)$}'' means a parameter range $0.4 \leq \zetaII \leq 0.9$, $0.9 \leq \omegaI \leq 1.1$, and $0.7 \leq \omegaII \leq 1.1$.
    \item Ultimately, we need to choose a set of parameters for use in practical settings. One approach could be to adopt the LFA-predicted best set of parameters. However, this optimum is sometimes situated off to one side and close to the boundary of the near-optimal region. Instead, we have chosen smoother parameters that better represent the average over all good parameters, i.e., at most 5\% above optimal; as such, we define the overall optimum as the centroid of the region $\bigl\{\rho(\zetaII,\omegaI,\omegaII)\leq \rho_\star^{1/1.05}\bigr\}$. In \cref{fig:optimal2,fig:optimal3}, this optimum is illustrated by the ball $\ball$, whose coordinates are given at the top of each figure.
\end{itemize}
\cref{app:lfa} provides further details on the LFA optimisation algorithm that adaptively explores the parameter space, prioritising areas where $\rho$ is small. More sophisticated approaches include the robust optimisation methods recently developed by Brown \textit{et al} \cite{doi:10.1137/19M1308669}; these methods specifically target multigrid optimisation problems, including the possibility of incorporating derivative and subgradient relations specific to LFA, and could lead to substantial speedups compared to the simple search algorithm used in this work.

We conclude this section by referring back to \cref{algo:Q} and noting that, for the remainder of this paper and its numerical experiments, we apply the optimal smoothing parameters as defined by the $\ball$ coordinates in \cref{fig:optimal2,fig:optimal3}. In particular, the left (resp., middle) column of \cref{fig:optimal2,fig:optimal3} defines the parameters $(\zetaII^0,\omegaI^0,\omegaII^0)$ for steady-state standard-form (resp., stress-form) Stokes problems, while the right column defines $(\zetaII^\infty,\omegaI^\infty,\omegaII^\infty)$.

\section{Numerical Experiments}
\label{sec:experiment}

In this section, we examine multigrid performance on a variety of Stokes problems in two and three dimensions and across multiple grid sizes, for a variety of polynomial degrees $\polydeg$ and boundary condition types. The metric used to assess multigrid efficiency concerns the convergence rate of the corresponding multigrid V-cycle left-preconditioned GMRES solver. Specifically, for each test case, we randomly generate a vector $b$ and solve $Ax = b$ with zero initial guess, $x_0 \equiv 0$. Each iteration $k$ of GMRES yields an improved solution $x_k$ and corresponding residual norm $r_k$; typically, $r_k$ is a nearly-constant fraction of $r_{k-1}$ so that $r_k \approx \rho^k\,r_{0}$, where $\rho$ is the convergence rate.\footnote{The components of $b$ are i.i.d.~randomly drawn from the uniform distribution on $[-1,1]$ so that the test problem is likely to contain modes for which convergence is the slowest. The residual norm is the same as that used by left-preconditioned GMRES, i.e., $r_k := \|V\!Ax_k - Vb\|_2$. On a log-linear plot, the data points $(k, r_k)$ are well approximated by a straight line of negative slope; simple linear regression is used to compute the corresponding convergence rate $\rho$. Our numerical tests typically have 5--10 such data points, except when convergence is so fast that just a few data points are collected before reaching the limits of double-precision arithmetic. {Transient behaviour is sometimes observed over the first 1--2 iterations, after which the preconditioned residual norm closely tracks that of the true residual, each with steady geometric decay.}} Our primary performance metric is then defined by
\begin{equation} \label{eq:eta} \eta := \frac{\log 0.1}{\log \rho} . \end{equation}
Here, $\eta$ represents the expected number of iterations needed to reduce the residual by a factor of $10$; e.g., $k\eta$ represents the number of iterations to reduce it by $10^k$. We prefer to use $\eta$ instead of $\rho$ as the former is more intuitively useful \cite{inferno}; e.g., a two-fold reduction in $\eta$ coincides with a two-fold reduction in iteration count, and so, all else being equal, two times faster.

For each prototype Stokes problem, we consider polynomial degrees $\polydeg \in \{1,2,3\}$ and test on grid sizes ranging from $4 \times 4$ up to $1024 \times 1024$ in 2D, and $4 \times 4 \times 4$ up to $256 \times 256 \times 256$ in 3D. We also test against three types of boundary conditions: (i) periodic, in which both $\vu$ and $p$ are periodic on $\partial \Omega$; (ii) Dirichlet, in which $\vu$ is prescribed on $\partial \Omega$ and $p$ is free; and (iii) stress, in which $\vu$ and $p$ are coupled via the specification of $\vsigma \cdot \vn$ on $\partial \Omega$. We test all three kinds as they can stress the numerics in rather different ways: (i) periodic boundary conditions are arguably the easiest to handle, and focus on the bulk-domain multigrid smoothing behaviour; (ii) velocity Dirichlet boundary conditions are enforced only weakly, via penalty operators that have very different spectral and smoothing characteristics compared to the bulk domain; and (iii) stress boundary conditions, which are not penalised at all, yet in some sense are a blend of Dirichlet boundary conditions on pressure and Neumann boundary conditions on the velocity field.

Our numerical results are presented through graphs of the form shown in \cref{fig:A}. Each column corresponds to a specific spatial dimension and polynomial degree and contains three sequences of markers, one for each boundary condition type: symbols $\bcP$, $\bcD$, and $\bcN$ denote periodic, Dirichlet, and stress boundary conditions, respectively. From left-to-right, the markers in each sequence plot the experimentally-determined multigrid convergence speed $\eta$ as measured on a grid of size $n \times n\,(\times\,n)$, where $n = 4, 8, 16, 32, \ldots$, in that order.

\subsection{Single-phase, steady-state, standard-form}
\label{sec:A}

\newcommand\orbscale{1}
\newcommand{\legend}{For each configuration (2D or 3D, polynomial degree $\polydeg$, and boundary condition type), the corresponding sequence of symbols plot the numerically determined speed $\eta$ on grid sizes $n \times n\,(\times\,n)$ where $n = 4,8,16,32,\ldots$, from left-to-right. Periodic, Dirichlet, and stress boundary conditions are denoted by $\bcP$, $\bcD$, and $\bcN$, respectively.}
\newcommand{\legendamr}{For each configuration (2D or 3D, polynomial degree $\polydeg$, and boundary condition type), the corresponding sequence of symbols plot the numerically determined speed $\eta$ on meshes with a base level grid size of $n \times n\,(\times\,n)$ where $n = 4,8,16,32,\ldots$, from left-to-right. Dirichlet and stress boundary conditions are denoted by $\bcD$ and $\bcN$, respectively.}

\begin{figure}%
\centering%
\includegraphics[scale=\orbscale]{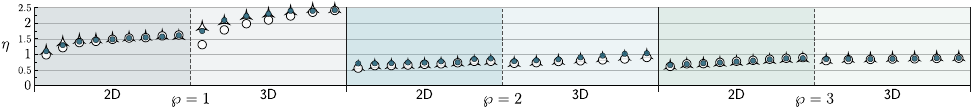}%
\caption{Multigrid solver performance for the single-phase steady-state standard-form Stokes problem considered in \cref{sec:A}. \legend}%
\label{fig:A}%
\end{figure}

We begin with perhaps the simplest case, that of a single-phase steady-state standard-form Stokes problem on the unit rectangular domain $\Omega = (0,1)^d$ with $\mu \equiv 1$, using uniform Cartesian grids. \Cref{fig:A} compiles the corresponding numerical results. For each combination of spatial dimension, polynomial degree, and boundary condition, we see that $\eta$ essentially plateaus, i.e., the multigrid iteration count required to reduce the residual by a fixed factor remains bounded as the grid is refined, as should be the case for an effective multigrid method. In particular, for polynomial degrees $\polydeg = 2$ and $\polydeg = 3$, a $10^{10}$-fold reduction in the residual requires 5 to 9 multigrid iterations---these speeds match those of classical geometric multigrid methods applied to Poisson problems. On the other hand, the case of $\polydeg = 1$ yields noticeably larger iteration counts, e.g., with $\eta$ approaching $2.5$ in 3D. This outcome is observed in most of our numerical results: for the same residual reduction threshold, the case of $\polydeg = 1$ requires 2--3 times as many outer solver iterations compared to the cases of $\polydeg \in \{2,3\}$. In this regard, $\polydeg = 1$ is a kind of outlier; on the other hand, the corresponding discrete solution is low-order (indeed, it only yields first-order pressure fields) and thus may not require many digits of accuracy, thereby allowing for a looser convergence criterion. Note also that in the context of high-order LDG methods, $\polydeg = 1$ is perhaps only of peripheral interest; further comments are given in the concluding remarks.

\subsection{Single-phase, steady-state, stress-form}
\label{sec:B}

\begin{figure}%
\centering%
\includegraphics[scale=\orbscale]{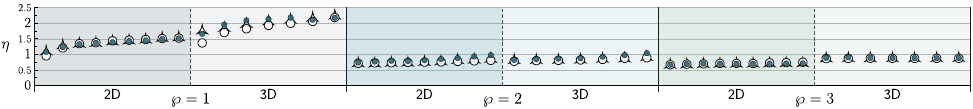}%
\caption{Multigrid solver performance for the single-phase steady-state stress-form Stokes problem considered in \cref{sec:B}. \legend}%
\label{fig:B}%
\end{figure}

We next test multigrid performance when considering steady-state Stokes equations in stress form, i.e., when $\gamma = 1$. We consider again a single-phase problem on $\Omega = (0,1)^d$ with $\mu \equiv 1$ and uniform Cartesian grids. \Cref{fig:B} compiles the corresponding numerical results. Compared to the standard-form case, we observe essentially the same multigrid efficiency, though with a mild speedup in some configurations; e.g., for $\polydeg = 3$ in 2D, $\eta$ plateaus to around $0.7$ (resp., $0.9$) for the stress-form (resp., standard-form) cases, see the second-from-right column of \cref{fig:B} (resp., \cref{fig:A}).

\subsection{Adaptive mesh refinement}
\label{sec:C}

\begin{figure}%
\centering%
\raisebox{-0.5\height}{\includegraphics[height=1.5in]{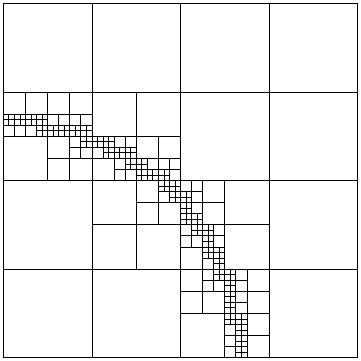}} \qquad \raisebox{-0.5\height}{\includegraphics[height=1.75in]{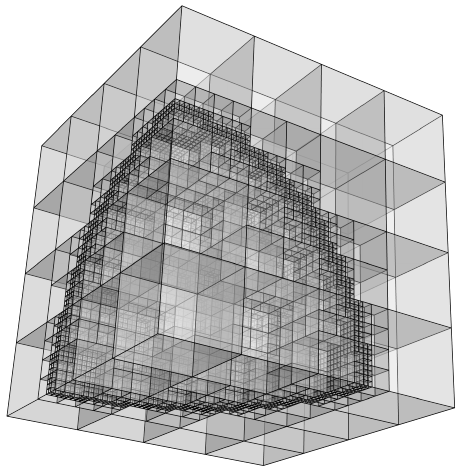}}%
\caption{Examples of the nonuniform quadtree (left) and octree (right) used in the adaptive mesh refinement problem in \cref{sec:C}. Each example corresponds to a base level grid size of $n = 4$, refined up to four levels around a circular arc or spherical shell of the same radius.}%
\label{fig:amr}%
\end{figure}

\begin{figure}%
\centering%
\includegraphics[scale=\orbscale]{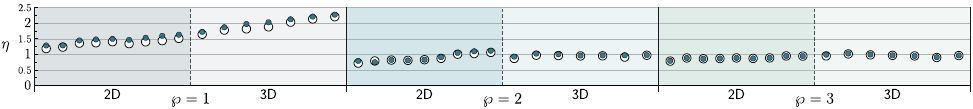}%
\caption{Multigrid solver performance for the quadtree/octree nonuniform mesh problem considered in \cref{sec:C}. \legendamr}%
\label{fig:C}%
\end{figure}

Our third test problem is similar to the first (in \cref{sec:A}), except that we now apply it to a nonuniform adaptively refined mesh. The mesh has up to four levels of nested refinement around a circular arc (in 2D) or spherical shell of the same radius (in 3D), as illustrated in \cref{fig:amr}. Owing to the rapid rate of refinement, this problem has considerably different stencil and sparsity structure. Compiled in \cref{fig:C} are the corresponding experimentally-determined multigrid solver speeds.\footnote{The largest such problem corresponds to the 3D case wherein the base level $256 \times 256 \times 256$ mesh is refined by four more levels, to a total of $\approx 88$ million elements; for $\polydeg = 3$, this corresponds to a velocity/pressure solution vector of $\approx 19$ billion degrees of freedom.} Examining \cref{fig:C}, we observe similar multigrid speeds as to prior test problems, though with a modest $\approx 25\%$ increase in $\eta$ that can be attributed to the more complex mesh geometry, stencils, and sparsity structure.

\subsection{Multiphase, steady-state, stress-form}
\label{sec:D}

\begin{figure}%
\centering%
\includegraphics[scale=\orbscale]{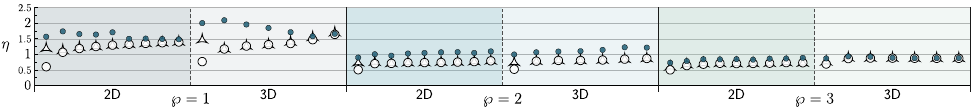}\\%
\includegraphics[scale=\orbscale]{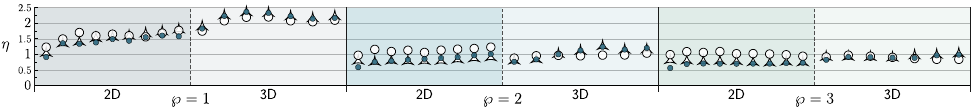}\\
\caption{Multigrid solver performance for the multiphase steady-state stress-form Stokes problems considered in \cref{sec:D}; the top (resp., bottom) row corresponds to an interior phase with viscosity a factor $10^4$ (resp., $10^{-4}$) that of the exterior phase. \legend}%
\label{fig:D}%
\end{figure}

We now consider a more challenging case of a high-contrast multiphase Stokes problem, wherein the viscosity $\mu$ exhibits a jump several orders in magnitude across an embedded interface. Let $\smash{\Omega = (0,1)^d}$ be divided into an interior phase, $\smash{\Omega_0 = (\tfrac14,\tfrac34)^d}$, and an exterior phase, $\smash{\Omega_1 = \Omega \setminus \overline{\Omega_0}}$, each with viscosity $\mu_i$. Note that in this multiphase example, there is an implicit specification of interfacial jump conditions, on both velocity and stress, on $\Gamma = \partial \Omega_0 \cap \partial \Omega_1$. Two representative viscosity ratios are examined, $\smash{\mu_0/\mu_1 = 10^4}$ and $\smash{\mu_0/\mu_1 = 10^{-4}}$; \cref{fig:D} compiles the corresponding numerical results. Compared to the prior test problems, we see a greater dependence on the applied domain boundary conditions. For example, when the interior phase is a factor $10^4$ more viscous than the exterior, we see from the top row of \cref{fig:D} that Dirichlet boundary conditions generally lead to larger $\eta$ compared to the other boundary conditions; this is most pronounced with $\polydeg = 1$ and to a lesser extent on $\polydeg = 2$, though with negligible difference for $\polydeg = 3$. Conversely, when the exterior phase is a factor $10^4$ more viscous than the interior phase, we see from the bottom row of \cref{fig:D} that sometimes stress boundary conditions lead to mildly slower convergence. Notwithstanding these effects, we observe fast multigrid convergence that remains stable under grid refinement; this favourable outcome is attributed in part to the use of viscosity-upwinded numerical fluxes, which bias the interfacial LDG discretisation to the mutual benefits of multigrid speed and high-order accuracy \cite{fluxx} as well as to the balancing process described in \cref{sec:balance}. Viscosity ratios of $10^{\pm2}$, $10^{\pm6}$, and $10^{\pm 8}$ have also been tested, with essentially identical results to those shown in \cref{fig:D}.

\subsection{Single-phase, unsteady, standard-form}
\label{sec:E}

So far we have focused on steady-state Stokes problems; our last set of test cases examine unsteady Stokes problems. As discussed in \cref{sec:unsteady}, the amended system includes a $\smash{\tfrac{\rho}{\dt}}$-weighted mass matrix whose strength, relative to the viscous operator, determines the nature of the unsteady Stokes problem: in one extreme, the viscous operator dominates and we expect to observe multigrid behaviour akin to the steady-state case; in the other extreme, the viscous operator is negligible so that we are in effect solving a flux-form Poisson problem with an effective ellipticity coefficient equal to $\smash{\tfrac{\dt}{\rho}}$. In between, the smoother parameters $(\zetaII,\omegaI,\omegaII)$ are blended; recall that the blending amount may actually change across the domain and across the multigrid hierarchy, since it depends on the local grid-dependent relative strengths of the $\rho$-weighted and $\mu$-weighted operators.

\begin{figure}%
\centering%
\includegraphics[scale=\orbscale]{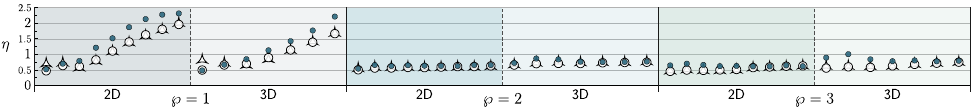}\\%
\includegraphics[scale=\orbscale]{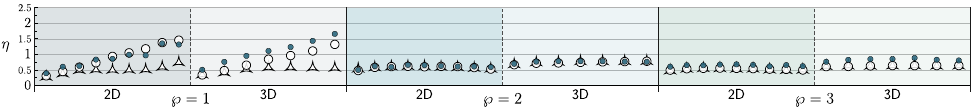}\\
\caption{Multigrid solver performance for the single-phase unsteady standard-form Stokes problems considered in \cref{sec:E}; the top (resp., bottom) row corresponds to an effective Reynolds number of $10^2$ (resp., $10^4$). \legend}%
\label{fig:E}%
\end{figure}

In the specific examples of this subsection, we consider two single-phase, unsteady, standard-form Stokes problems on $\Omega = (0,1)^d$: one with $\mu = 10^{-2}$ and the other with $\mu = 10^{-4}$. In both cases, we set $\rho = 1$ and define $\dt = 0.1h$ where $h$ is the element size of the primary-level uniform Cartesian grid; this choice of $\dt$ represents a typical scenario of applying the unsteady Stokes equations in a temporal integration method having unit-order CFL number. These two problems correspond to an effective Reynolds number $\Rey = \rho\,U\!L/\mu$ of $\Rey = 10^2$ and $\Rey = 10^4$, respectively: the former (with $\mu = 10^{-2}$) represents a viscous-dominated case, but where the time-derivative term nevertheless influences solver characteristics; the latter (with $\mu = 10^{-4}$) represents a case wherein the time-derivative term dominates, though with progressively weaker influence as the primary mesh is refined. With this setup, \cref{fig:E} compiles the numerical results. For polynomial degrees $\polydeg \in \{2,3\}$, we observe relatively stable and fast multigrid solver performance for all kinds of boundary condition; in particular, $\eta$ is often close to $0.5$, i.e., just $5$ multigrid GMRES iterations are needed to reduce the residual by ten orders of magnitude. Meanwhile, the case of $\polydeg = 1$ exhibits more nuanced behaviour: for the high-Reynolds case with periodic boundary conditions, the multigrid solver performs well, reaching speeds of $\eta \approx 0.5$; however, for all other cases we observe an upwards trend in $\eta$. This is largely attributed to the fact that, as the grid is refined, eventually we obtain a perturbed unsteady Stokes problem, for which we have seen in prior examples (see \cref{fig:A,fig:B,fig:C,fig:D}) that $\eta \approx 2$ is typical.

\subsection{Multiphase, unsteady, stress-form}
\label{sec:F}

\begin{figure}%
\centering%
\includegraphics[scale=\orbscale]{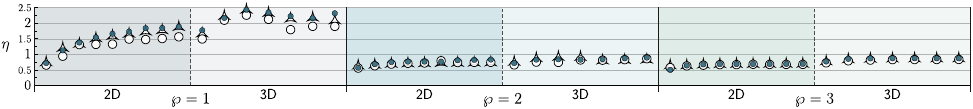}\\%
\includegraphics[scale=\orbscale]{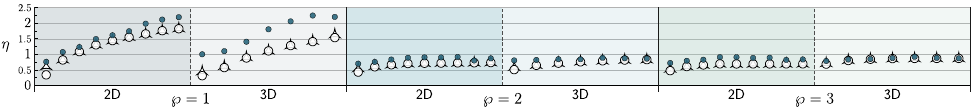}\\
\caption{Multigrid solver performance for the multiphase unsteady stress-form Stokes problems considered in \cref{sec:F}; the top (resp., bottom) row corresponds to a gas droplet surrounded by liquid (resp., a liquid droplet surrounded by gas). \legend}%
\label{fig:F}%
\end{figure}

Our last two examples are in some regards the most challenging, as they involve density and viscosity coefficients that jump by 3 to 4 orders of magnitude across an embedded interface, thereby examining the ability to handle high-contrast coefficients in addition to the grid-dependent relative strengths of the time derivative and viscous operator terms. We consider two scenarios: a liquid bubble surrounded by gas, and a gas bubble surrounded by liquid, with $\rho_{\text{liquid}} = 1$, $\mu_{\text{liquid}} = 1$, $\rho_{\text{gas}} = 0.001$, and $\mu_{\text{gas}} = 0.0002$ (approximately accurate values in CGS units for ambient air and 100 cSt silicone oil). Here, the ``bubble'' refers to an interior rectangular phase, $\smash{\Omega_0 = (\tfrac14,\tfrac34)^d}$, while the exterior phase is $\smash{\Omega_1 = \Omega \setminus \overline{\Omega_0}}$, with $\Omega = (0,1)^d$. Note this problem also examines the multigrid solver's ability to handle velocity and stress interfacial jump conditions across the embedded liquid-gas interface. As before, we set $\dt = 0.1h$ to represent, say, solving the incompressible Navier-Stokes equations via an explicit treatment of the advection term and implicit treatment of the viscous term. With this setup, \cref{fig:F} compiles the corresponding numerical results. As in other test problems, the low-order case of $\polydeg = 1$ is a kind of outlier; excluding it, we observe multigrid speeds that attain those of classical geometric multigrid methods: essentially grid-independent convergence of around 5 to 9 solver iterations to reduce the residual by ten orders of magnitude.

\section{Concluding Remarks}
\label{sec:conclusion}

The two main LDG approaches for Stokes problems are that of equal-degree (whereby $\vu_h \in V_\polydeg^d$ and $p_h \in V_\polydeg$) and mixed-degree (whereby $\vu_h \in V_\polydeg^d$ and $p_h \in V_{\polydeg-1}$). Both yield optimal order accuracy for the velocity field, i.e., order $\polydeg + 1$ in the maximum norm; and, except for idealised cases, both yield order $\polydeg$ accuracy in the maximum norm for the pressure field.\footnote{It is possible for the equal-degree approach to yield order $\polydeg+1$ in the maximum norm, for both velocity and pressure, but only in some idealised cases such as periodic boundary conditions on uniform grids using one-sided fluxes.} However, the equal-degree approach often leads to significant pressure errors near nontrivial boundaries and interfaces, whereas the mixed-degree approach yields considerably more accurate and cleaner pressure fields, see \cref{sec:comparison}. But, and perhaps superseding aspects of accuracy, is the complication of pressure penalty stabilisation. In general, pressure stabilisation is required when using the equal-degree approach, but there is no known efficient way to implement it for stress boundary conditions nor interfacial jump conditions, see \cref{sec:comparison}. In contrast, the mixed-degree approach does not require  pressure stabilisation, thereby eliminating many subtle implementation details.

Besides the development of a mixed-degree LDG framework for multiphase Stokes problems, in this paper we focused mainly on multigrid solvers. In particular, we developed a multigrid smoother that is identical to multicoloured element-wise block Gauss-Seidel, except the block diagonal inverses $A_{ii}^{-1}$ are replaced by approximate inverses $Q_i$. In turn, this inverse is computed through a weighted least squares approximation of the form $A Q \approx \mathbb{I}$, controlled by four tunable parameters, $\zetaI, \zetaII, \omegaI$, and $\omegaII$. The $\zeta$ parameters control the weights in the least squares problem, while the $\omega$ parameters allow the velocity and pressure variables to be damped/overrelaxed by different amounts. Intuitively, their combined purpose is to simultaneously and maximally damp the high frequency error components of the velocity and pressure variables, in part by balancing the relative strengths of the stress divergence and velocity divergence operators. In addition, the smoother applies a balancing algorithm to treat two critical aspects of Stokes problems: (i) the mesh-dependent, unequal scaling of the viscous, pressure gradient, velocity divergence, and density-weighted temporal derivative operators; and (ii) for high-contrast multiphase problems, viscosity and/or density coefficients varying over multiple orders in magnitude. In effect, this balancing process augments the weighted least squares formulation so that the aforementioned smoother parameters $(\zetaI, \zetaII, \omegaI, \omegaII)$ can be applied independent of physical parameters, element sizes, mesh connectivity, etc.

The smoother parameters can be tuned in various ways. Fortunately, $\zetaI$ can be quickly eliminated from consideration, essentially by choosing it arbitrarily large, see \cref{sec:lfa}. For the remaining three parameters, we used an optimisation  approach based on two-grid local Fourier analysis (LFA). The outcome is an optimised set of parameters that depend solely on the spatial dimension $d$, the polynomial degree $\polydeg$, and the form of the Stokes problem (whether it is steady-state standard-form, steady-state stress-form, or unsteady vanishing-viscosity); these parameters are listed in \cref{fig:optimal2,fig:optimal3}. Extensive numerical tests---over a variety of steady and unsteady Stokes problems, in 2D and 3D, across multiple grid sizes, on polynomial degree $\polydeg \in \{1,2,3\}$, and with different boundary conditions---demonstrated that the resulting multigrid solver matches the speed of {Poisson-style} geometric multigrid methods, e.g., 5 to 10 iterations to reduce the residual by a factor of $10^{10}$ (putting aside the low-order outlier of $\polydeg = 1$).

\begin{table}%
\centering%
\sffamily\footnotesize%
\begin{tabular}{lcc|ccccc}
& & & \multicolumn{3}{c}{Relative cost $mn^2/n^3$} \\
Stokes problem & $d$ & $\ell$ & $\polydeg = 1$ & $\polydeg = 2$ & $\polydeg = 3$ \\
\midrule
\multirow{2}{*}{Standard form}
& 2D & 4 & 1.89 & 1.36 & 1.20 \\
& 3D & 6 & 1.72 & 1.20 & 1.08 \\
\midrule
\multirow{2}{*}{Stress form}
& 2D & 6 & 2.33 & 1.55 & 1.29 \\
& 3D & 12 & 2.44 & 1.40 & 1.16 \\
\midrule
\end{tabular}%
\caption{Relative cost of computing $Q_i$ vs.~$A_{ii}^{-1}$, for canonical Stokes problems on uniform $d$-dimensional Cartesian grids.}%
\label{tab:relative}
\end{table}

Regarding the smoother's computational efficiency, this was designed from the outset to be the same or comparable to the equivalent block Gauss-Seidel method. {In particular, a single iteration is directly comparable to standard block Gauss-Seidel in per-iteration cost, memory footprint, and scalability.} Meanwhile, the dominant component of building $Q_i$ is a least squares problem of size $m \times n$, where $n = d(\polydeg + 1)^d + \polydeg^d$ is the elemental block size and $m = n + d \ell$, where $\ell$ equals the number of nonzero off-diagonal blocks (e.g., a standard-form Stokes problem in 3D has a 7-element stencil, so $\ell = 6$). To leading order, the corresponding solution requires $\mathcal{O}(m n^2)$ floating point operations, whereas computing $A_{ii}^{-1}$ requires $\mathcal{O}(n^3)$ operations; \cref{tab:relative} provides a rough indication of their relative cost by tabulating the ratio $(mn^2)/n^3$ for canonical uniform-grid Stokes problems. In practice, the cost of precomputing $Q$ is almost negligible---e.g., to reduce the residual by a factor of $10^{10}$, say, the time spent precomputing $Q$ is no more than 10\% of the overall multigrid solver time. Additionally, note that each block of $Q$ can be precomputed independently of every other block, making its construction trivially parallelisable and scalable.

As observed in our numerical results, the case of $\polydeg = 1$ is an outlier in the sense it requires 2--3 times more iterations compared to the multigrid solvers for $\polydeg > 1$. However, this comparison is perhaps incomplete as it overlooks some application details; e.g., low-order solutions are unlikely to require many digits of accuracy, so the same solver tolerances need not apply. %
Some numerical experiments (not shown here) suggest that the $\polydeg = 1$ smoother parameters could be tweaked for an additional 10--20\% speedup; in other words, the predictions of two-grid LFA appear to be somewhat suboptimal for $\polydeg = 1$. To that end, a three-grid LFA method may have better predictive power in determining the absolute-best smoother parameters. Alternatively, generalising the approach used here to an overlapping-block smoother method, similar to Vanka-type smoothers, offers a more promising avenue for building effective multigrid solvers in the $\polydeg = 1$ case; we leave this exploration open for future study. %

Finally, we mention some extensions and future research possibilities. We considered here an LDG modal polynomial basis. This choice mainly impacted the design and construction of the smoother, where it was used to easily select a few rows per off-diagonal block of $A_{:,i}$ (i.e., the rows corresponding to the constant-mode lowest-degree modal coefficient of the stress divergence operator). If a nodal basis was preferred, one could make a temporary nodal-to-modal change of basis, build the block-diagonal smoother $Q$, and then map it back to the nodal basis. An alternative approach might be to preserve the nodal basis, but select some other subset of rows from the off-diagonal blocks---e.g., rows corresponding to collocated nodes on neighbouring elements---though this approach would likely require a new set of optimised smoother parameters. On a similar note, we expect these ideas could be adapted to other LDG settings, e.g., triangular/tetrahedral meshes, total-degree polynomial bases, etc., and perhaps other kinds of DG discretisations besides LDG. Meanwhile, we focused here on V-cycles with 3 pre-smoothing and 3 post-smoothing steps. Other possibilities include V- and W-cycles with differing pre- and post-smoother counts, F-cycles, etc. As suggested by Brown \textit{et al} \cite{doi:10.1137/19M1308669}, even faster multigrid solvers could be developed by applying one set of smoother parameters in the pre-smoother and a different set of parameters in the post-smoother.

\section*{Acknowledgements}

This work was supported in part by the U.S. Department of Energy, Office of Science, Office of Advanced Scientific Computing Research's Applied Mathematics Competitive Portfolios program under Contract No.~DE-AC02-05CH11231, and by a U.S. Department of Energy, Office of Science, Early Career Research Program award. Some computations used resources of the National Energy Research Scientific Computing Center (NERSC), a U.S. Department of Energy, Office of Science User Facility operated at Lawrence Berkeley National Laboratory.

\appendix

\section{Local Discontinuous Galerkin Framework}%
\label{app:ldg}%
We derive here a mixed-degree LDG framework for multiphase variable-viscosity Stokes problems. The construction closely follows that of prior work \cite{flame}, except herein the pressure field has polynomial degree one less than the velocity field. At the same time, the framework is partly based on the LDG schemes initially developed by Cockburn \textit{et al} \cite{CockburnKanschatSchotzauSchwab2002}, extending them to multiphase variable-viscosity problems, incorporating interfacial jump conditions, viscosity-upwinded fluxes, and additional types of domain boundary conditions.

For simplicity of presentation, we consider meshes $\mathcal E$ arising from Cartesian grids, quadtrees, or octrees. To establish some notation, \textit{intraphase} mesh faces are those shared by two elements in the same phase; \textit{interphase} mesh faces are those shared by two elements in differing phases (and thus located on $\Gamma_{ij}$ for some $i,j$); and \textit{boundary} mesh faces are those situated on $\partial \Omega$. Each face has an associated unit normal vector: interphase faces adopt the same unit normal as their corresponding interface $\Gamma_{ij}$; boundary faces adopt the natural outwards-pointing domain boundary unit normal; while intraphase faces are all positively oriented in the same way, e.g., $\vn = \hat{\mathbf x}$ on every vertical intraphase face. The notation $\jump{\cdot}$ denotes the jump of a quantity across an interface and/or mesh face: $\jump{u} := u^- -  u^+$ where $u^\pm(x) = \lim_{\epsilon \to 0^+} u(x \pm \epsilon \vn)$ evaluates the trace of $u$ on either side. In addition, for an element $E \in \mathcal E$, let $\vn_E$ denote the outwards pointing normal on the boundary of the element; define $\chi(E)$ to be the phase of $E$ (so that $E \subseteq \Omega_{\chi(E)}$); and for any discontinuous piecewise polynomial function $u$, let $u|_E$ denote the restriction of $u$ to $E$ (so that the quantity $u|_E$ is well defined on $\partial E$). 

In the first step of the LDG construction, a discrete approximation to $\nabla \vu$ is defined. Given $\vu \in V_\polydeg^d$, we define ${\bm \tau} \in V_\polydeg^{d \times d}$ to be the discrete gradient of $\vu$ via a ``strong-weak'' form, so that
\begin{equation} \label{eq:tau} \int_E {\bm \tau} : {\bm \omega} = \int_E \bm \omega : \nabla \vu + \int_{\partial E} ({\hat \vu}_{\chi(E)} - \vu|_E) \cdot \bm \tau|_E \cdot \bm n_E \end{equation}
holds for every element $E$ and every test function $\bm \omega \in V_\polydeg^{d \times d}$. Here, the numerical flux $\hat{\vu}_\chi$ is defined by
\begin{equation} \label{eq:uhat} \hat{\vu}_\chi := \begin{cases}
\vu^- & \text{on any intraphase face,} \\
\lambda\,\vu^- + (1 - \lambda)\,(\vu^+ + \vg_{\chi j}) & \text{on $\Gamma_{\chi j}$ if $\chi < j$,} \\
\lambda\,(\vu^- - \vg_{i\chi}) + (1 - \lambda)\,\vu^+ & \text{on $\Gamma_{i\chi}$ if $i < \chi$,} \\
\vu^- & \text{on $\Gamma_N$,} \\
\vg_\partial & \text{on $\Gamma_D$}.
\end{cases} \end{equation}
On intraphase faces, we have chosen a minimal dissipation LDG scheme \cite{CockburnKanschatSchotzauSchwab2002}, whereby the numerical flux is purely one-sided, being sourced from the ``left'' element, i.e., $\hat\vu = \vu^-$. On interphase faces, the numerical flux is multivalued so as to implement the interfacial jump condition $\jump{\vu} = \vg_{ij}$; whenever the weak formulation of \cref{eq:tau} ``reaches across'' an interface $\Gamma_{ij}$ to evaluate the trace of $\vu$ on the other side, the result is offset by the jump data $\vg_{ij}$ so as to appropriately account for the intended discontinuity in the solution. The blending of the trace values on either side of the interface is controlled by the convex weighting parameter $\lambda \in [0,1]$; its specification allows for the implementation of viscosity-upwinded fluxes \cite{fluxx}.

A routine calculation (see \cite{flame} for an example) shows that, upon summing \cref{eq:tau} over every mesh element, the variational form is equivalent to the statement that
\[ \tau_{ij} = G_j u_i + J_{\vg,ij}, \]
where $\smash{G = (G_1, \ldots, G_d) : V_\polydeg \to V_\polydeg^d}$ is a discrete gradient operator implementing the numerical fluxes of \cref{eq:uhat}, $\tau_{ij} \in V_\polydeg$ denotes the $(i,j)$th component of $\bm \tau$, $u_i$ denotes the $i$th component of $\vu$, and $J_\vg \in V_\polydeg^{d \times d}$ implements the lifting of any Dirichlet-like source data $\vg_\partial$, $\vg_{ij}$. (The precise definitions of $G$ and $J_\vg$ are summarised in \cref{sec:ldgdef}.)

In the second step of the LDG construction, a discrete approximation of $\vsigma = \mu (\bm \tau + \gamma\,\bm \tau^\trans) - p\,\mathbb I$ is defined. Two aspects need attention: (i) the viscosity $\mu$ may vary spatially (and may possibly be a tensor, though for simplicity of presentation we do not consider that here); and (ii) the pressure field $p$ is of polynomial degree one less than the velocity field. Both aspects can be concisely treated by performing an $L^2$ projection of $\mu (\bm \tau + \gamma\,\bm \tau^\trans) - p\,\mathbb I$ onto $V_\polydeg^{d \times d}$, resulting in the definition of $\vsigma$ as
\begin{equation} \label{eq:sigmaresult} \sigma_{ij} := M^{-1} M_\mu \bigl(\tau_{ij} + \gamma\,\tau_{ji} \bigr) - \delta_{ij}\,{\mathbb P}_{V_\polydeg}(p) = M^{-1} M_\mu \bigl( G_j u_i + \gamma\,G_i u_j + J_{\vg,ij} + \gamma\,J_{\vg,ji} \bigr) - \delta_{ij}\,{\mathbb P}_{V_\polydeg}(p). \end{equation}
Here, $M_\mu$ is a $\mu$-weighted mass matrix relative to $V_\polydeg$, while ${\mathbb P}_{V_\polydeg}(p)$ is the projection of $p \in V_{\polydeg-1}$ onto $V_\polydeg$ and simply corresponds to injecting $p$ into the higher-degree polynomial space.

In the third step, we consider the weak formulation for computing the discrete divergence of $\vsigma$. This proceeds similarly to the first step, except numerical fluxes are taken in the opposite direction. (For simplicity of presentation, the following numerical flux for $\vsigma$ is matrix-valued, however only the normal component is used.) Define $\bm \xi \in V_\polydeg^d$ as the discrete divergence of $\vsigma$ so that
\begin{equation} \label{eq:xidef} \int_E \bm \xi \cdot \bm \omega = - \int_E \vsigma : \nabla \bm \omega + \int_{\partial E} \bm \omega|_E \cdot \hat{\vsigma}_{\chi(E)} \cdot \vn_E \end{equation}
holds for every element $E$ and every test function $\bm \omega \in V_\polydeg^d$. Here, the numerical flux $\hat \vsigma$ is defined by
\begin{equation} \label{eq:sigmahat} \hat{\vsigma}_\chi := \begin{cases}
\vsigma^+ & \text{on any intraphase face,} \\
(1 - \lambda)\,\vsigma^- + \lambda\,(\vsigma^+ + \vh_{\chi j} \otimes \vn) & \text{on $\Gamma_{\chi j}$ if $\chi < j$,} \\
(1 - \lambda)\,(\vsigma^- - \vh_{i \chi} \otimes \vn) + \lambda\,\vsigma^+ & \text{on $\Gamma_{i \chi}$ if $i < \chi$,} \\
\vh_\partial \otimes \vn & \text{on $\Gamma_N$,} \\
\vsigma^- & \text{on } \Gamma_D.
\end{cases} \end{equation}
Analogous to the numerical flux $\hat{\vu}$, the interfacial jump condition $\jump{\vsigma \cdot \vn} = \vh_{ij}$ is taken into account in \cref{eq:sigmahat} such that, whenever an element reaches across the interface, the neighbouring element's trace is offset by $\vh_{ij}$ so as to appropriately account for the intended discontinuity in $\vsigma \cdot \vn$. Summing \cref{eq:xidef} over every mesh element, a routine calculation (see \cite{flame} for an example) shows that the variational form leads to
\begin{equation} \label{eq:nueqn} \xi_i = -\sum_{j=1}^d M^{-1} G_j^\trans M \sigma_{ij} + J_{\vh,i}, \end{equation}
where $\xi_i \in V_\polydeg$ is the $i$th component of $\bm \xi$ and $J_{\mathbf h} \in V_\polydeg^d$ (defined in \cref{sec:ldgdef}) implements the lifting of any Neumann-like source data $\vh_{\partial}$, $\vh_{ij}$. Note that $\xi_i$ computes the discrete divergence of $\vsigma$ via the adjoint of the discrete gradient operator $G$; the adjoint arises because the direction and weighting of the numerical fluxes for $\hat\vsigma$ in \cref{eq:sigmahat} are opposite to those of $\hat\vu$ in \cref{eq:uhat}. To complete this step, we (temporarily) set the negative divergence of $\vsigma$ equal to the $L^2$ projection of the Stokes momentum equation source data $\vf$; combining \cref{eq:sigmaresult} and \cref{eq:nueqn} we therefore obtain, for every component $i$,
\begin{equation} \label{eq:divsigma} \sum_{j=1}^d M^{-1} G_j^\trans M_\mu \bigl( G_j u_i + \gamma\,G_i u_j + J_{\vg,ij} + \gamma\,J_{\vg,ji} \bigr) - M^{-1} G_i^\trans M {\mathbb P}_{V_\polydeg}(p) - J_{\vh,i} = {\mathbb P}_{V_\polydeg^d}({\vf}_i). \end{equation}
One may identify the first term as a discrete approximation to $-\nabla \cdot \bigl(\mu (\nabla \vu + \gamma\,\nabla \vu^\trans) \bigr)$ and the second as a mixed-degree discrete gradient of pressure.

We now turn to the LDG formulation for the Stokes divergence constraint. Given $\vu \in V_\polydeg^d$, we define $\varphi \in V_{\polydeg-1}$ as a discrete divergence of $\vu$ via a strong-weak form, so that
\begin{equation} \label{eq:divu} \int_E \varphi\,\omega = \int_E \omega \, \nabla \cdot \vu + \int_{\partial E} \omega \, (\hat\vu_{\chi(E)} - \vu|_E) \cdot \vn_E \end{equation}
holds for every element $E$ and every test function $\omega \in V_{\polydeg-1}$. Here, the same numerical flux $\hat \vu$ is used, i.e., \cref{eq:uhat}, but only the normal component is seen by the weak formulation in \cref{eq:divu}. Summing \cref{eq:divu} over the entire mesh, it follows that
\[ \varphi = \sum_{i = 1}^d \widetilde{G}_i u_i + J_{\vg \cdot \vn}, \]
where $\widetilde{G} = (\widetilde{G}_1, \cdots, \widetilde{G}_d) : V_\polydeg \to V_{\polydeg-1}^d$ is a mixed-degree discrete gradient operator and $J_{\vg \cdot \vn} \in V_{\polydeg-1}$ implements the lifting of the normal component of any Dirichlet-like source data $\vg_\partial$, $\vg_{ij}$. Note that $\widetilde{G}$ maps the full-degree space to the lower-degree space of the pressure field; in fact, $\widetilde{G}$ is identical to $G$ except for an $L^2$ projection, i.e., $\widetilde{G}_i = {\mathbb P}_{V_{\polydeg-1}} G_i$. 
To complete this step, the negative discrete divergence of $\bm u$ is set equal to the $L^2$ projection of the Stokes divergence equation source data $f$, thereby obtaining
\begin{equation} \label{eq:divuresult} -\sum_{i = 1}^d \widetilde{G}_i u_i - J_{\vg \cdot \vn} = {\mathbb P}_{V_{\polydeg-1}}(f). \end{equation}

Finally, in the last step of the LDG construction, penalty stabilisation terms are included so as to ensure the well-posedness of the discrete problem. Solely targeting the velocity field, these terms weakly enforce Dirichlet boundary conditions and interfacial jump conditions; they may also be used to weakly enforce continuity of the velocity field between neighbouring elements of the same phase. Let $\tau_{\vu,\partial}$, $\tau_{\vu,ij}$, and $\tau_{\vu,\circ}$ denote the corresponding penalty parameters for the boundary, interphase, and intraphase faces, respectively. Let $E_{\vu,\vg} : V_\polydeg^d \to V_\polydeg^d$ be the affine operator such that
\[ \int_\Omega E_{\vu,\vg}(\vu) \cdot \bm \omega = \int_{\Gamma_D} \tau_{\vu,\partial} \, (\vu^- - \vg_\partial) \cdot \bm \omega^- + \sum_{i < j} \int_{\Gamma_{ij}} \tau_{\vu,ij}\, \bigl( \jump{\vu} - \vg_{ij}  \bigr) \cdot \jump{\bm \omega} + \int_{\Gamma_\circ} \tau_{\vu,\circ}\, \jump{\vu} \cdot \jump{\bm \omega} \]
holds for every test function $\bm \omega \in V_\polydeg^d$. (Here, $\Gamma_\circ$ denotes the union of all intraphase faces.) Subject to a suitable specification of penalty parameter values, $E_{\vu,\vg}$ is added to the left-hand side of \cref{eq:divsigma}. Multiplying the resulting equations by $M$, multiplying \cref{eq:divuresult} by $\bar M$, and then rearranging the resulting system of equations so that the entire influence of the source data is moved to the right-hand side, we then obtain the governing equations \cref{eq:components} as described in the main article.

\subsection{Operator definitions}
\label{sec:ldgdef}

For completeness, we define here the various operators and source terms omitted from the above derivation.

\begin{itemize}
    \item The mass matrix $M$ relative to the full-degree space $V_\polydeg$ is defined such that $u^\trans M v = \int_\Omega u\,v$ for every $u, v \in V_\polydeg$; similarly for the mass matrix $\bar M$ relative to the space $V_{\polydeg-1}$. Meanwhile, the $\mu$-weighted mass matrix $M_\mu$ is defined such that $u^\trans M_\mu v = \sum_i \int_{\Omega_i} u\,\mu_i\,v$, for all $u, v \in V_\polydeg$.
    \item The discrete gradient operator $G : V_\polydeg \to V_\polydeg^d$ is defined by $G := \bg + L$, where the broken gradient operator $\bg : V_\polydeg \to V_\polydeg^d$ and lifting operator $L : V_\polydeg \to V_\polydeg^d$ are defined such that 
    \[ \int_\Omega \bigl(\bg u\bigr) \cdot \bm\omega = \sum_{E \in \mathcal E} \int_E (\nabla u) \cdot \bm \omega, \quad \text{and} \]
    \[ \int_\Omega (Lu) \cdot \bm \omega = \int_{\Gamma_\circ} (u^+ - u^-)\, {\bm \omega}^+ \cdot \vn - \int_{\Gamma_D} u^- {\bm \omega}^- \cdot \vn + \sum_{i < j} \int_{\Gamma_{ij}} (1 - \lambda)\,(u^+ - u^-)\,{\bm \omega}^- \cdot \vn + \lambda\,(u^+ - u^-)\,{\bm \omega}^+ \cdot \vn \]
    hold for every $u \in V_\polydeg$ and $\bm\omega \in V_\polydeg^d$.
    \item The mixed-degree discrete gradient operator $\widetilde{G} : V_\polydeg \to V_{\polydeg-1}^d$ is defined as the $L^2$ projection of $G$ onto $V_{\polydeg-1}^d$, i.e., so that
    \[ \int_\Omega \bigl(\widetilde{G} u \bigr) \cdot \bm \omega = \int_\Omega (Gu) \cdot \bm \omega \]
    holds for every $u \in V_\polydeg$ and $\bm\omega \in V_{\polydeg-1}^d$. Written differently, note that $u^\trans \widetilde{G}_i^\trans \bar M \omega = u^\trans G_i^\trans M {\mathbb P}_{V_\polydeg} (\omega)$ for each component $i$ and every $u \in V_\polydeg$ and $\omega \in V_{\polydeg-1}$. The latter form shows that the pressure gradient term in \cref{eq:divsigma} can also be rewritten as $M^{-1} \widetilde{G}_i^\trans \bar M p$.
    \item The lifted Dirichlet-like source data $J_\vg \in V_\polydeg^{d \times d}$ is defined such that 
    \[ \int_\Omega J_\vg : \bm \omega = \int_{\Gamma_D} \vg_\partial \cdot \bm \omega^- \cdot \vn + \sum_{i < j} \int_{\Gamma_{ij}} (1 - \lambda)\, \vg_{ij} \cdot \bm \omega^- \cdot \vn + \lambda\, \vg_{ij} \cdot \bm \omega^+ \cdot \vn \]
    holds for every $\bm \omega \in V_\polydeg^{d \times d}$, while $J_{\vg \cdot \vn} \in V_{\polydeg-1}$ is defined such that
    \[ \int_\Omega J_{\vg \cdot \vn}\,\omega = \int_{\Gamma_D} \omega^-\,\vg_\partial \cdot \vn + \sum_{i < j} \int_{\Gamma_{ij}} \bigl((1 - \lambda)\,\omega^- + \lambda\,\omega^+ \bigr)\,\vg_{ij} \cdot \vn\]
    holds for every $\omega \in V_{\polydeg-1}$. Meanwhile, the lifted Neumann-like source data $J_\vh \in V_\polydeg^d$ is defined such that
    \[ \int_\Omega J_\vh \cdot \bm \omega = \int_{\Gamma_N} \vh_\partial \cdot \bm \omega^- + \sum_{i < j} \int_{\Gamma_{ij}} \lambda\,\vh_{ij} \cdot \bm \omega^-  + (1-\lambda)\,\vh_{ij} \cdot \bm \omega^+ \]
    holds for every $\bm \omega \in V_\polydeg^d$. 
    
\end{itemize}

\subsection{Viscosity-upwinded fluxes}

For multiphase Stokes problems, we apply viscosity-upwinded numerical fluxes, developed for elliptic interface problems in \cite{fluxx} and extended to multiphase Stokes problems in \cite{flame}. Roughly speaking, in the limit of arbitrarily-large viscosity ratios, this biasing corresponds to the multiphase Stokes problem separating into two distinct problems: one for the highly viscous phase (which effectively ``sees'' a stress-like boundary condition at the interface), and one for the other phase (which effectively ``sees'' a Dirichlet-like boundary condition). As explored in \cite{fluxx}, biasing the numerical fluxes in the same way is a key step in achieving efficient multigrid solvers as well as high-order accuracy. The viscosity-upwinded strategy corresponds to defining $\lambda$ in the numerical flux functions $\hat\vu$ and $\hat\vsigma$ in \cref{eq:uhat} and \cref{eq:sigmahat} as follows:
\begin{equation} \label{eq:viscosityupwind} \lambda = \begin{cases} 0 & \text{if } \mu^- < \mu^+, \\ 0.5 & \text{if } \mu^- = \mu^+, \\ 1 & \text{if } \mu^- > \mu^+. \end{cases} \end{equation}
Accordingly, $\hat\vu$ is biased to the more viscous phase and $\hat\vsigma$ to the less viscous phase, with interfacial jump data $\vg_{ij}$ and $\vh_{ij}$ incorporated appropriately.

\subsection{Well-posedness}

Assuming a suitable specification of penalty parameters (see below), the mixed-degree LDG framework is observed to be well-posed, inf-sup stable, and optimal-order accurate. The formulation may be viewed as a generalisation of the LDG schemes developed by Cockburn \textit{et al.} \cite{CockburnKanschatSchotzauSchwab2002}, for which optimal-order accuracy and inf-sup stability have been proven for single-phase, constant-coefficient model problems; see also \cite{CockburnKanschatSchotzau2003,CockburnKanschatSchotzau2004}. Although we have not formally proven these properties for the generalised setting, extensive numerical experiments demonstrate that they hold, including for smoothly varying viscosity coefficients as well as high-contrast multiphase problems. We leave a formal proof for future work.

\subsection{Penalty stabilisation parameters}

In general, strictly positive penalty parameters $\tau_{\vu,\partial}$, $\tau_{\vu,ij}$, and $\tau_{\vu,\circ}$ are sufficient to ensure the well-posedness of the final discretised Stokes problem, e.g., to ensure it has the expected kernel of the continuum Stokes operator. %
However, strict positivity is not a necessary condition; e.g., on regular Cartesian grids, without any embedded interfaces, it is possible to entirely disable intraphase penalisation and still maintain stability. Further discussion on particular parameter choices and appropriate scaling with respect to viscosity, mesh size, polynomial degree, etc., are given in \cite{flame}; here, we specify the values used in the numerical experiments of this work:
\begin{itemize}
    \item On Dirichlet boundary faces, $\tau_{\vu,\partial} = C_\polydeg \mu^- h^{-1}$, where $C_\polydeg$ is a constant depending on the polynomial degree, $\mu^-$ is the local viscosity, and $h$ the local element size. For the numerical experiments in this work, $C_1 = 1$ while $C_2 = C_3 = 16$; the case of $\polydeg = 1$ benefits from somewhat weaker Dirichlet boundary penalisation.
    \item The AMR test problems require some amount of intraphase stabilisation, where it is set to $\tau_{\vu,\circ} = C_\polydeg h^{-1}$. In all other test problems, intraphase penalisation is unnecessary and is disabled, $\tau_{\vu,\circ} = 0$.
    \item For the multiphase test problems, the interphase penalty parameter is set to $\tau_{\vu,ij} = 16\,\min(\mu^-,\mu^+)\, h^{-1}$, where $\mu^\pm$ is the local viscosity of the two elements on either side of the interfacial mesh face. We note that it is important to choose the minimum of the two viscosities, as discussed further in \cite{fluxx}.
\end{itemize}
These values have been chosen based on the dual goals of obtaining high-order accuracy as well as good multigrid performance, and follow typical values employed in previous work \cite{fluxx,flame,inferno}.

\section{Optimising Smoother Parameters}%
\label{app:lfa}%
We describe here an algorithm to optimise for the multigrid smoother parameters $(\zetaII, \omegaI, \omegaII)$. The approach has two main components: (i) a search algorithm that adaptively explores the parameter space; and (ii) a method to evaluate the corresponding multigrid convergence rate $\rho$ via two-grid local Fourier analysis.

\begin{figure}%
\centering%
\begin{minipage}{0.48\textwidth}%
\begin{algorithm}[H]
    \caption{\sffamily\small Exploring the $(\zetaII, \omegaI, \omegaII)$ parameter space.}
    \begin{algorithmic}[1]
        \State Define initial seed point, $\mathcal{Z} := \bigl\{ (1, 1, 1) \bigr\}$.
        \State Define initial step size, $\delta := (0.25, 0.25, 0.25)$.
        \State Let $\varsigma$ track the minimum observed value of $\rho$ so that at any instant, $\varsigma \equiv \min_{z \in \mathcal{Z}} \rho(z)$.
        \While{exploring} \label{alg:outerloop}
            \State $\rho_\star := \begin{cases} \varsigma^\alpha & \text{if } \varsigma < 1, \\ +\infty & \text{otherwise.} \end{cases}$ \label{alg:rhostar}
            \State $\mathcal{Z}_\star := \bigl\{ z \in \mathcal{Z} \mid \rho(z) \leq \rho_\star\bigr\}$.
            \State $n_i := \frac{1}{\delta_i}\displaystyle\Bigl(\max_{z \in \mathcal{Z}_\star} z_i - \min_{z \in \mathcal{Z}_\star} z_i \Bigr), \forall i \in \{1, 2, 3\}$.
            \If{$\min_i n_i \geq n_\text{target}$}
                \State End search and \textbf{return} $\mathcal{Z}$.
            \EndIf
            \State Refine step size: $\delta_j := \tfrac12 \delta_j$, where $j = \argmin_i n_i$.
            \State Mark every point in $\mathcal{Z}$ as unvisited.
            \While{$\mathcal{Z}$ has unvisited points} \label{alg:innerloop}
                \State $z := \displaystyle\argmin_{\text{unvisited}\,\mathsf{z} \in \mathcal{Z}} \rho(\mathsf{z})$.
                \State Mark $z$ as visited.
                \If{$\varsigma < 1$ \textbf{and} $\rho(z) > \varsigma^\alpha$}
                    \State \textbf{break}
                \EndIf
                \For{$y \in \displaystyle\prod_{i = 1,2,3} \{z_i - \delta_i, z_i, z_i + \delta_i\}$}
                    \If{$y \notin \mathcal{Z}$}
                        \State Insert $y$ into $\mathcal{Z}$ (marked unvisited).
                    \EndIf
                \EndFor
            \EndWhile
        \EndWhile
\end{algorithmic}%
\label{algo:search}%
\end{algorithm}%
\end{minipage}%
\hfill%
\begin{minipage}{0.48\textwidth}%
\begin{algorithm}[H]
    \caption{\sffamily\small Evaluation of $\rho(\zetaII, \omegaI, \omegaII)$.}
    \begin{algorithmic}[1]
        \State Let $\Theta_n$ denote the set $\bigl\{ \frac{2\pi}{n} (i - \tfrac12) \bigr\}_{i=1}^n$.
        \State $\rho_\text{coarse} := \displaystyle\max_{\bm \theta \in (\Theta_6)^d} \rho( \bm\theta; \zetaII, \omegaI, \omegaII \bigr)$.
        \If{$\rho_\text{coarse} > 1$}
            \State \textbf{return} $\rho_\text{coarse}$.
        \EndIf
        \If{$\varsigma < 1$ \textbf{and} $\rho_\text{coarse} > \varsigma^\alpha$}
            \State \textbf{return} $\rho_\text{coarse}$.
        \EndIf
        \State $\rho_\text{fine} := \displaystyle\max_{\bm \theta \in (\Theta_{12})^d} \rho( \bm\theta; \zetaII, \omegaI, \omegaII \bigr)$.
        \State \textbf{return} $\max(\rho_\text{coarse}, \rho_\text{fine})$.
\end{algorithmic}%
\label{algo:rho}%
\end{algorithm}%
\end{minipage}%
\end{figure}

Suppose for the moment that we have implemented a function $\rho : \R^3 \to \R$ that evaluates the multigrid convergence rate $\rho(z)$ for the parameters $z = (\zetaII, \omegaI, \omegaII)$. Treating $\rho$ as a black-box function, \Cref{algo:search} outlines the search algorithm that iteratively and adaptively evaluates $\rho$ on a point cloud in parameter space. Given an (initial or interim) point cloud $\mathcal{Z}$ arising from a step size $\delta$, the method selects from $\mathcal{Z}$ the set of near-optimal points, $\mathcal{Z}_\star$, i.e., those for which the corresponding convergence rate $\rho$ is close to the minimum thus far. Next, the bounding box of $\mathcal{Z}_\star$ is computed so as to determine the effective sampling resolution relative to the step size $\delta$. If the sampling resolution across all three parameters is sufficiently high, the search process terminates. Otherwise, the coordinate axis with the lowest resolution is chosen to undergo refinement, so that $\delta$ is halved in that direction. Then, for each point in $\mathcal{Z}$---visited in order of increasing $\rho$ value---a patch of neighbouring points is constructed via a $3 \times 3 \times 3$ Cartesian grid of cell size $\delta$, and these neighbouring points are added to the search space. The sampling procedure prioritises areas where $\rho$ is expected to be the smallest; if the creation of a neighbouring point leads to the reduction of $\min_{z \in \mathcal Z} \rho(z)$, that point will be visited next. Conversely, if upon visiting an existing point whose $\rho$ value is now considered too large, all remaining unvisited points have the same property, so the inner loop (line \ref{alg:innerloop}) search can terminate. The outer loop (line \ref{alg:outerloop}) progressively applies finer resolutions; once the target resolution is reached, the overall search terminates, with the result being a point cloud $\mathcal{Z}$ for which the majority of points focus in on the region where $\rho$ is close to the minimum-observed $\rho$.

To specify a few more details in this search process, a sample $z$ is considered close to optimal if $\rho(z) \leq \bigl( \min_{\mathsf{z} \in \mathcal{Z}} \rho(\mathsf{z}) \bigr)^\alpha$, where $\alpha < 1$ is a user-chosen parameter; for the results in this work, we have used $\alpha = 1/1.1$, corresponding to the statement that the smoother parameters at $z$ yield at most 10\% more multigrid iterations (for a fixed residual reduction threshold) compared to the best found so far. The sampling process is initialised at the point $(\zetaII,\omegaI,\omegaII) = (1,1,1)$ with an initial step size of $0.25$. (Sometimes the initial seed yields a  $\rho$ larger than unity; when this occurs, the near-optimal threshold  $\rho_\star$ is set to $+\infty$ (line \ref{alg:rhostar}), thereby switching to an ``explore everywhere'' heuristic.) As the search proceeds, it may be the case that the typical range of $\zetaII$, say, ends up being larger than that of $\omegaI$, say; the convergence criterion roughly requires that the bounding box of the set of near-optimal points has at least $n_\text{target}$ points per dimension; we have used $n_\text{target} = 16$, sufficient to produce the figures of \cref{fig:optimal2,fig:optimal3}. Our implementation of \cref{algo:search} memoises/caches the results of the (expensive) black-box evaluation of $\rho$ in a simple dictionary/map data structure; meanwhile, the inner loop over unvisited points ordered by $\rho$ can be efficiently implemented via a simple priority queue.

We now turn to the evaluation of $\rho(z)$ for a given $z = (\zetaII, \omegaI, \omegaII)$. We define $\rho$ via the results of two-grid local Fourier analysis (LFA) \cite{brandt,mgbook1,mgbook2,doi:10.1137/19M1308669}. In essence, this procedure analyses the performance of a two-level multigrid method as applied to an infinite uniform Cartesian grid spanning all of $\R^d$, together with a function space spanned by functions of the form $\exp (\iota \bm \theta \cdot \bm x/h)$, where $\iota^2 = -1$. LFA analyses often proceed by recognising that a constant coefficient PDE operator is diagonalised by the Fourier modes $\exp (\iota \bm \theta \cdot \bm x/h)$, thereby yielding ``symbols'' of the PDE operator $A$ as well as symbols of the multigrid pre-smoothing, post-smoothing, restriction, interpolation, and coarse grid operators. These symbols combine to form a symbol $E(\bm \theta)$ for the error propagation operator for a two-level method---in the present setting, $E(\bm \theta)$ is a complex-valued square matrix of size proportional to the \#dof on each element---from which the convergence rate, for a specific frequency $\bm \theta$, is defined as the spectral radius of $E(\bm \theta)$. The maximum spectral radius over an appropriate frequency space is then used to define the LFA-predicted convergence rate $\rho$ for the overall multigrid method. For details, see, e.g., \cite{mgbook1,mgbook2,doi:10.1137/19M1308669}.

In this work, instead of symbolically building and manipulating the various symbols, we found it simpler to mimic two-grid LFA by invoking the multigrid DG implementation on a virtual infinite-extent Cartesian grid of cell size $h$. Given a frequency $\bm\theta \in [-\pi,\pi)^d$, we define a virtual grid function $\bm x \mapsto \sigma \exp(\iota\,\bm \theta \cdot \bm x / h)$, where $\bm x$ denotes the lower-left coordinates of each grid cell. Here, $\sigma$ denotes a single set of coefficients of size $d (\polydeg + 1)^d + \polydeg^d$ (one per velocity and pressure modal coefficient, respectively); the coefficients $\sigma$ are the same for every grid cell, but otherwise remain undefined. For a given set of smoother parameters $(\zetaII, \omegaI, \omegaII)$, the multigrid smoother is computed following the methods of \cref{sec:mg}. Then, a two-level multigrid method is applied to solve $A_h x_h = b_h$, as follows (for simplicity of presentation we consider the 2D case, the 3D case being an obvious extension):
\begin{enumerate}
    \item Define an initial grid function $x_h$ on the virtual infinite-extent fine grid, $x_h^{i,j} = \sigma \exp\bigl(\iota\, \bm \theta \cdot (i,j) \bigr)$ on grid cell $(i,j) \in \mathbb Z^d$. Define a right-hand side of $b_h \equiv 0$.
    \item Apply the presmoother $\nu_1$ times, using red-black ordering. The result of smoothing is a new grid function such that $x_h^{i,j} = (X_{i,j} \sigma) \exp\bigl(\iota\,\bm \theta \cdot (i,j) \bigr)$, where $X_{i,j}$ is cell-dependent complex-valued matrix. Crucially, owing to the nature of red-black smoothing, the new grid function repeats every $2 \times 2$ patch of grid cells, in the sense that $X_{i,j} = X_{i \bmod 2, j \bmod 2}$ for every $(i,j) \in \mathbb Z^2$.
    \item Compute the residual on the fine grid, $r_h = b_h - A_h x_h$; the result is another grid function which repeats every $2 \times 2$ grid cells.
    \item Restrict the residual to the virtual infinite-extent coarse grid; crucially, the restriction process agglomerates every patch of $2 \times 2$ cells into a single element, which means the restricted residual $b_{2h}$ is the same on every coarse grid cell, in the sense that $b_{2h}^{i,j} = (B \sigma) \exp\bigl(\iota\,\bm \theta \cdot (i,j)\bigr)$ where $B$ is a cell-independent complex-valued matrix and $(i,j) \in (2\mathbb Z)^2$.
    \item Solve the coarse-grid problem, exactly: using the symbol of $A_{2h}$, compute $x_{2h}$ such that $A_{2h} x_{2h} = b_{2h}$. The result is a grid function such that $x_{2h}^{i,j} = (C \sigma) \exp\bigl(\iota\,\bm \theta \cdot (i,j)\bigr)$ for every coarse grid cell $(i,j) \in (2\mathbb Z)^2$; in fact, $C = A_{2h}(\bm \theta)^{-1} B$, where $A_{2h}(\bm \theta)$ is the symbol of the coarse grid operator and $A_{2h}(\bm \theta)^{-1}$ can be computed, e.g., via complex-valued SVD.
    \item Interpolate the coarse-grid solution and apply it as a correction to the fine grid function $x_h$. Then apply the postsmoother $\nu_2$ times. The updated fine grid function is again of the form $x_h^{i,j} = (X_{i,j} \sigma) \exp\bigl(\iota \bm \theta \cdot (i,j) \bigr)$, where $X_{i,j}$ is a cell-dependent complex-valued matrix; as before, these matrices repeat: $X_{i,j} = X_{i \bmod 2, j \bmod 2}$ for every $(i,j) \in \mathbb Z^2$, so only four such matrices need computing.
    \item At this point, we have in fact computed the two-grid error propagator of the two-level multigrid method, represented as a linear function applied to the unknown coefficients $\sigma$. Finally, we determine the slowest convergence rate over all possible coefficients, i.e., compute the spectral radius of the error propagator. Using the matrices of the previous step, we therefore define the two-grid multigrid convergence rate, for the input frequency $\bm \theta$ and smoother parameters $(\zetaII,\omegaI,\omegaII)$, via
    \[ \rho(\bm \theta; \zetaII,\omegaI,\omegaII) := \max_{(i,j) \in \{0,1\}^2} \rho\bigl(X_{i,j} \bigr). \]
\end{enumerate}
The above method evaluates the two-grid convergence factor $\rho$ for a concrete, numerically-defined frequency $\bm \theta$. Its implementation is nearly identical to the regular multigrid method, except for three main aspects: (i) on each element, state variables are complex-valued matrices instead of the usual real-valued coefficient vector; (ii) stencil operations appropriately account for the grid function complex-exponential factors, e.g., the $\begin{bmatrix} -1 & 2 & -1 \end{bmatrix}$ stencil of a 1D finite difference Laplacian is replaced by $\begin{bmatrix} -e^{-\iota \theta} & 2 & -e^{\iota\theta} \end{bmatrix}$; and (iii) instead of applying the multigrid algorithm to every cell, because of the repeating nature of the computations, we need only apply it to a $2 \times 2\,(\times\,2)$ primary patch of cells. Finally, we find the slowest convergence rate across all possible frequencies to define the final two-grid LFA convergence factor:\footnote{In typical LFA approaches, it is common to consider blocked forms of the error propagation symbol, blocked according to a primary low-frequency space $[-\tfrac\pi2,\tfrac\pi2)^d$, together with multiple $\pi$-shifted copies thereof, and then maximise over the low-frequency space \cite{doi:10.1137/19M1308669}. For the alternative approach presented in this paper, it is appropriate to maximise over the entire frequency space.}
\begin{equation} \label{eq:rhothetaparams} \rho(\zetaII,\omegaI,\omegaII) := \sup_{\bm \theta \in [-\pi,\pi)^d} \rho(\bm \theta; \zetaII,\omegaI,\omegaII). \end{equation}
In this work, we have used a simple maximisation approach, outlined in \cref{algo:rho}, that works in tandem with the parameter search of \cref{algo:search}, as follows. As a preliminary coarse approximation, the above two-grid LFA method is used to compute $\rho(\bm \theta; \zetaII,\omegaI,\omegaII)$ on a $6 \times 6\,(\times\,6)$ grid in frequency space. If the maximum-found $\rho$ is too large compared to the best smoother parameters found so far (located somewhere else), then it is in some sense futile to perform a more accurate calculation. If, on the other hand, the coarse approximation results in a sufficiently small $\rho$, then we apply a higher-resolution evaluation on a $12 \times 12\,(\times\,12)$ grid in frequency space. The net effect is the evaluation of \cref{eq:rhothetaparams} on semi-uniform grid in frequency space having $6+12$ points per dimension; in this work it was found that the coarse resolution of $6$ does well to quickly rule out bad smoother parameters, while the fine resolution of $6+12$ yields a predicted convergence factor $\rho$ sufficiently close to the true maximum, at least for the overall purpose of predicting a good set of smoother parameters as discussed in \cref{sec:lfa}.

\section{A Brief Comparison with Classical Smoothers}%
\label{app:briefcomparison}%
We consider here an illustrative two-dimensional model problem to compare multigrid performance across a variety of classical smoothers. The model problem consists of a single-phase, constant-viscosity, steady-state, standard-form Stokes problem, with periodic boundary conditions, on an $n \times n$ uniform Cartesian grid, and $\polydeg = 2$, i.e., biquadratic velocity field and bilinear pressure field. The solver consists of a multigrid V-cycle left-preconditioned GMRES method, with $\nu = 3$ pre- and post-smoothing steps. We consider three kinds of classical smoothers: (i) block Jacobi; (ii) multicoloured block Gauss-Seidel using identical pre- and post-smoothing colour sequences (forward/forward); and (iii) multicoloured block Gauss-Seidel using a reversed post-smoothing colour sequence (forward/reverse). For each of these, we consider undamped and optimally damped versions.\footnote{A simple parameter sweep was used to select an approximately optimal damping factor $\omega$ that maximises the convergence rate.} \Cref{tab:briefcomparison} reports the number of outer GMRES iterations required to achieve a residual reduction factor of $10^8$ (using a zero initial guess and a randomly generated right-hand side); blank entries denote failure to converge within 50 iterations. Overall, block Jacobi is ineffective, while multicoloured block Gauss-Seidel degrades as the mesh is refined. In contrast, the $Q$ smoother developed in this work maintains fast and nearly grid-independent convergence. Other model problem configurations, such as three-dimensional cases, yield similar findings.

\begin{table}%
\centering%
\sffamily\footnotesize%
\begin{tabular}{p{6.5cm}|c|cccccccc}
& & \multicolumn{7}{c}{Grid size $n$} \\
Multigrid smoother & Damping & 4 & 8 & 16 & 32 & 64 & 128 & 256 & 512 \\
\midrule
\multirow{2}{*}{Block Jacobi}
& none              & -- & -- & -- & -- & -- & -- & -- & -- \\
& $\omega = 0.55$   & 13 & 16 & 18 & 21 & 24 & 32 & 40 & 47 \\
\midrule
\multirow{2}{=}{Block Gauss-Seidel (forward/forward)}
& none              & -- & -- & -- & -- & -- & -- & -- & -- \\
& $\omega = 0.7$    & 7 & 8 & 9 & 10 & 10 & 13 & 16 & 21 \\
\midrule
\multirow{2}{=}{Block Gauss-Seidel (forward/reverse)}
& none              & 50 & -- & -- & -- & -- & -- & -- & -- \\
& $\omega = 0.7$    & 8 & 8 & 10 & 10 & 12 & 13 & 15 & 17 \\
\midrule
Least-squares balanced block-diagonal $Q$ (this work)
& n/a               & 5 & 5 & 5 & 5 & 6 & 6 & 6 & 6 \\
\midrule
\end{tabular}%
\caption{Multigrid solver iteration counts for the 2D model problem considered in \cref{app:briefcomparison}; blank entries denote failure to converge within 50 iterations.}%
\label{tab:briefcomparison}
\end{table}

\bibliographystyle{siamplain}
\bibliography{references}

\end{document}